%% file: hfk2.tex
\documentclass[11pt]{article}
\usepackage{amsmath, amsbsy, amscd, amsthm, amsfonts, array, epic, eepic, verbatim}
 \theoremstyle{definition}

 \numberwithin{equation}{section}

\def\ZZ{\mathbb{Z}}

\def\Sym{\operatorname{Sym}}

\begin{document}  

\title{COMPUTATIONS OF HEEGAARD-FLOER \\ KNOT HOMOLOGY}
\author{John A.~Baldwin and W.~D.~Gillam \\ Columbia University Department of Mathematics}

\date{September, 2006}

\maketitle

\begin{abstract} Using a combinatorial approach described in a recent paper of Manolescu, Ozsv\'ath, and Sarkar we compute the Heegaard-Floer knot homology of all knots with at most 12 crossings as well as the $\tau$ invariant for knots through 11 crossings.  We review the basic construction of \cite{MOS}, giving two examples that can be worked out by hand, and explain some ideas we used to simplify the computation.  We conclude with a discussion of knot Floer homology for small knots, closely examining the Kinoshita-Teraska knot $KT_{2,1}$ and its Conway mutant.
\end{abstract}

\begin{section}{Introduction} \label{section:intro}
In \cite{OSz5}, Ozsv\'ath and Szab\'o introduced a topological invariant of closed, orientable three-manifolds in the form of a collection of abelian groups. To such a three-manifold $Y$, we associate the groups (we use $\ZZ_2$ coefficients below, so these will be $\ZZ_2$ vector spaces) $\widehat{HF}(Y,s)$, indexed by $Spin^c$ structures on $Y$. In \cite{OSz3}, they discovered that a null-homologous knot $K$ in a three-manifold $Y$ induces a (finite) filtration $$ \cdots \subseteq \widehat{CFK}^i(Y,K) \subseteq \widehat{CFK}^{i+1}(Y,K) \subseteq \cdots \subseteq \widehat{CF}(Y)$$ on the chain complex used to compute the Heegaard Floer homology groups of $Y$.\footnote{This was independently discovered by Jacob Rasmussen in his thesis \cite{R}.}  Moreover, the filtered chain homotopy type of this complex is an invariant of the knot $K \subset Y$. The homology of the successive quotients $$\widehat{HFK}_j(Y,K,i):=H_j( \widehat{CFK}^i(Y,K) / \widehat{CFK}^{i-1}(Y,K) )$$ (in each $Spin^c$ structure) is called the Heegaard-Floer knot homology of $K \subset Y$ and, until very recently, there was no algorithmic way to compute these knot Floer homology groups---computing boundary maps in the chain complex involves counting pseudo-holomorphic disks satisfying certain boundary conditions in $\Sym^g(\Sigma_g)$. 

In their remarkable paper \cite{MOS}, Ozsv\'ath, Manolescu, and Sarkar describe an algorithm for computing $\widehat{CFK}(S^3,K)$ for any knot $K \subset S^3$ (for a knot $K \subset S^3$ we will denote the corresponding knot Floer homology group simply by $\widehat{HFK}(K)$). The combinatorial description of knot Floer homology in \cite{MOS} begins with the observation that an arc presentation of a knot $K$ naturally gives rise to a genus one multi-pointed Heegaard diagram $D$ for the knot in $S^3$ for which the generators, gradings, filtration, and boundary maps in the associated filtered complex $C(D)$ are completely combinatorial. The number $n$ of arcs in this arc presentation will be equal to the number of $\alpha$ and $\beta$ curves in the Heegaard diagram mentioned above (hence $n$ can be made equal to the arc-index $\alpha(K)$ of $K$).  The filtered chain homotopy type of $C(D)$ is equal to $\widehat{CFK}(K) \otimes S^{\otimes(n-1)}$, where $S$ is the filtered complex over $\ZZ_2$ with trivial boundary operator given below. 

$$S^i_j = \left \{ \begin{array}{ll} \ZZ_2 & j = -1 {\rm \; and \;} i \geq -1 \\ \ZZ_2 & j=0 {\rm \; and \;} i \geq 0 \\ 0 & {\rm otherwise} \end{array} \right.$$ By computing the homology of successive quotients of $C(D)$ we can recover $\widehat{HFK}(K)$.

An arc presentation $D$ consists of an $n \times n$ grid with white and black dots placed inside squares of the grid in such a way that each row and each column of the grid contains exactly one white dot and one black dot (there is at most one dot in a square). This corresponds to an oriented knot (or link) $K$ as follows: in each column, connect the black and white dots by a vertical line oriented from the black dot to the white dot; in each row, connect the black and white dots by a horizontal line from the white dot to the black dot (always passing the horizontal strands under the vertical strands).  Then $C(D)$ is the free $\ZZ_2$-module on permutations of $\{1,\dots,n\}$.  We can identify a permutation $\sigma$ with its graph $\Gamma(\sigma)$ on the grid (indexing squares by their lower left corners; view the grid as being on a torus $T$ so the top and bottom lines are identified and the left and right sides are identified).  Label the circles on the torus determined by the vertical lines of the grid $\alpha_1,\dots,\alpha_n$ and the circles determined by the horizontal lines $\beta_1,\dots,\beta_n$.  Then the generators of $C(D)$ correspond to the points of $$[\alpha_1 \times \dots \times \alpha_n] \cap [\beta_1 \times \dots \times \beta_n] \subseteq \Sym^n T.$$  

We assign an Alexander grading to each permutation as follows.  For a point $p$ on the grid, let $w(p)$ denote the winding number of $K$ around $p$. To compute the winding number of a nice curve $C$ around a point $p$ in the plane, draw any ray $L$ from $p$ to $\infty$ and take the signed sum of intersection points $C \cap L$ (with the sign convention chosen so that a counterclockwise circle has winding number 1 around a point in the bounded region defined by the circle). Each black or white dot is in a square whose four corners are four points on the grid, giving us a total of $8n$ distinguised grid points $p_1,\dots,p_{8n}$ (there may be repeats).  Let $$a:= (1/8)\sum_{i=1}^{8n} w(p_i) - (n-1)/2$$ and declare the Alexander grading of a permutation $\sigma$ to be $$A(\sigma)=a+\sum_{p \in \Gamma(\sigma)} -w(p).$$  The filtration level $C^i(D) \subseteq C(D)$ is generated by permutations with Alexander grading at most $i$.  

Every permutation is also assigned a Maslov (homological) grading. To begin with, if $p$ is a point on the grid or in the interior of one of its squares, and $S$ is a square of the grid, let $S.p$ be 1 if $p \in S$ and $0$ otherwise. Let $\sigma_0$ denote the permutation whose graph consists of the lower left corners of the squares on the grid containing a white dot.  To compute the Maslov grading of a permutation $\sigma$, find the smallest $i_1 \in \{1,\dots,n\}$ for which $\sigma(i_1) \neq \sigma_0(i_1)$; there must be some $i_2$ so that $\sigma_0(i_2)=\sigma(i_1)$ so draw a horizontal line (always draw lines on the grid, forgetting that it is on a torus) from $(i_1,\sigma(i_1))$ to $(i_2,\sigma_0(i_2))$ and then a vertical line from this point to $(i_2,\sigma(i_2))$; then there must be some $i_3$ so that $\sigma_0(i_3)=\sigma(i_2)$ so draw the horizontal line from $(i_2,\sigma(i_2))$ to $(i_3,\sigma_0(i_3))$, and so on, continuing in this manner until a closed, oriented curve is produced.  If there is still some $i$ with $\sigma(i) \neq \sigma_0(i)$ for which $(i,\sigma(i))$ is not on this closed curve, take the smallest such $i$ and repeat the above process.  This eventually yields a collection of oriented closed curves $\gamma_{\sigma \sigma_0}$ on the grid, none of which wraps around the torus.  Thus $\gamma_{\sigma \sigma_0}$ is the oriented boundary of a formal sum $\sum_i a_i S_i$ of squares of the grid, and we declare the Maslov grading of $\sigma$ to be $$M(\sigma) = 1-n +\frac{1}{4} \sum_{i, \; p \in \Gamma(\sigma)} a_i S_i.p + \frac{1}{4} \sum_{i, \; p \in \Gamma(\sigma_0)} a_i S_i.p - 2 \sum_{i, \; {\rm white \; dots \;} p} a_i S_i.p.$$  Note that the mod $2$ Maslov grading only depends on the sign of the permutation $\sigma$.  

The boundary map $d:C_j(D) \to C_{j-1}(D)$ is defined as follows.  If $\sigma$ and $\sigma'$ differ at more than two places then neither appears as a boundary of the other.  If $\sigma$ and $\sigma'$ differ exactly on $i<j$, say with $\sigma(i)<\sigma'(i)$ then the circles $\alpha_i,\alpha_j,\beta_{\sigma(i)},\beta_{\sigma(j)}$ divide the torus into four rectangular connected components $R_1,R_2,R_3,R_4$ as in Figure~\ref{figure:rect}. 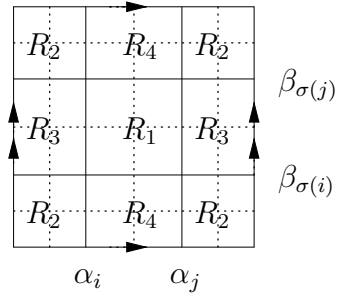
\begin{figure}[htbp] \begin{center}
\input{rect.eepic}
\vspace*{-2mm} \caption{Four rectangles used to compute boundaries}
\label{figure:rect}
\end{center}
\vspace*{-5mm}\end{figure}  Then $\sigma' \in d\sigma$ iff exactly one of $R_1,R_2$ contains no white dot, or point in the graph of $\sigma$ in its interior and $\sigma \in d \sigma'$ iff exactly one of $R_3,R_4$ contains no white dot, or point in the graph of $\sigma$ in its interior.  It is fairly straight-forward to check that $d^2=0$ and that this boundary map reduces Maslov grading by 1 and does not increase the Alexander grading.  Indeed, to compute the homology of successive quotients we simply disregard boundaries with lower Alexander grading.  This is equivalent to requiring the rectangles to be free of black dots as well.

\end{section}

\begin{section}{Methodology} \label{section:methodology}

We implemented the combinatorial description of $C(D)$ roughly as described above.  Then we used a fairly well-known algorithm (see Figure~\ref{figure:chain}) for determining generators for the homology of a complex $C$ over $\ZZ_2$. \begin{figure}[htbp]
\begin{center}
\input{chain.eepic}
\vspace*{-2mm}\caption{Reducing a chain complex over $\ZZ_2$}
\label{figure:chain}
\end{center}
\vspace*{-5mm}\end{figure}
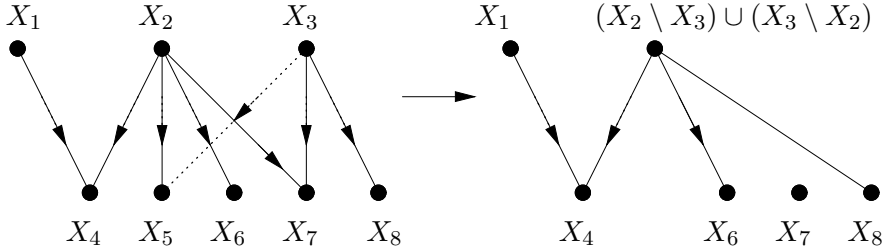  Choose a basis $B=\{x_1,\dots,x_m \}$ for $C$ and identify each element $\sum_{i=1}^n x_{j_i}$ ($j_1 \leq \dots \leq j_n$) of $C$ with the set $\{x_{j_1}, \dots, x_{j_n} \} \subseteq B$.  Start with the (directed) graph $\Gamma_0$ with vertex set $$\{ \{ x_1 \}, \dots, \{ x_m \} \}$$ and with an arrow from $\{ x_i \}$ to $\{ x_j \}$ exactly when $x_j \in dx_i$.  Now suppose we have a directed graph $\Gamma_n$ whose vertices are labelled with sets $X_i \subset B$ $(i \in I)$ in such a way that

\smallskip

$(*) \;\;\; \{ X_i : i \in I \}$ is linearly independent in $C$

\smallskip

\noindent and suppose there is an edge from $X_i$ to $X_j$ in this graph.  Then we produce another graph $\Gamma_{n+1}$ satisfying $(*)$ with two fewer vertices as follows.  The vertex set of $\Gamma_{n+1}$ will be $\{ Y_k : k \in I \setminus \{ i,j \} \}$ where $Y_k$ is the symmetric difference of $X_i$ and $X_k$ if $X_k$ is the start of an edge ending at $X_j$ and $Y_k=X_k$ otherwise.  In the former case there will be an edge from $Y_k$ to $Y_l$ in $\Gamma_{n+1}$ whenever exactly one of $X_k$ or $X_i$ is the start of an edge terminating at $X_l$ in $\Gamma_n$.  In the latter case there will be an edge from $Y_k$ to $Y_l$ exactly when there is an edge from $X_k$ to $X_l$ in $\Gamma_n$.   

Starting with $\Gamma_0$ we perform this algorithm repeatedly to obtain a sequence of graphs $\Gamma_0, \dots, \Gamma_N$ where $\Gamma_N$ has no edges.  Then the generators of the homology of $C$ are the vertices of $\Gamma_N$.  This is very fast when there are relatively few boundaries, as occurs in the situation above where there are $n!$ generators but at most $\begin{pmatrix} n \\ 2 \end{pmatrix}$ boundaries for each generator.  In practice, when $n=10$ each generator appears to have an average of approximately $7$ boundaries, so there are surprisingly few edges in the graph.  Because of the symmetry \cite{OSz3} $$\widehat{HFK}_j(S^3,K,i) =\widehat{HFK}_{j-2i}(S^3,K,-i)$$ it is enough to compute knot Floer homology in non-negative Alexander grading; again this cuts the computation down immensely since, of the $n!$ generators, only a small fraction have non-negative Alexander grading.  We do not compute Maslov grading until the graph is reduced as this computation is quite time consuming. 

Recall that to any (say, finitely generated) filtered chain complex $C$ one can associate a spectral sequence whose $E^1$ term is $E^1_{p,q} = H_{p+q}(C^p / C^{p-1} )$ converging to the homology of $C$.  Furthermore, the first differential $d^1:E^1_{p,q} \to E^1_{p-1,q}$ in this spectral sequence is identified with the connecting homomorphism in the LES in homology associated to the SES of complexes $$0 \to C^{p-1}/C^{p-2} \to C^p/C^{p-2} \to C^p / C^{p-1} \to 0.$$ When $C = \widehat{CFK}(K)$, we have $$E^1_{p,j-p} = \widehat{HFK}_j(K,p)$$ and the spectral sequence converges to the Heegaard-Floer homology of $S^3$ which is just a $\ZZ_2$ in homological grading $0$ so we may ask for the smallest integer $\tau$ such that the natural map $$H_0( \widehat{CFK}^{\tau}(K) ) \to H_0( \widehat{CF}(S^3) ) \cong \ZZ_2 $$ is non-zero.  This is the well-known concordance invariant defined by Ozsv\'ath and Szab\'o in \cite{OSz6}.  It turns out that $|\tau(K)| \leq g^*(K)$, where $g^*(K)$ is the smooth four-ball genus of $K$.

Fortunately the knot Floer homology of a knot with at most 11 crossings is small enough (i.e. is supported on at most 2 diagonals) that we can figure out $\tau$ just by computing the ranks of the $d^1$ differentials.  To compute the rank of $$d^1 : \widehat{HFK}_j(K,i) \to \widehat{HFK}_{j-1}(K,i-1)$$ we found explicit generators for $H_j( C^i(D) / C^{i-1}(D) )$ using the previously-described algorithm, then computed the part of their boundary in one lower Alexander grading to find cycles $$z_1, \dots, z_n  \in ( C^{i-1}(D) / C^{i-2}(D) )_{j-1}.$$  We then formed a new complex $B(D)$ from $C(D)$ by adding generators $g_1, \dots, g_n$ to $C(D)_{j-2}$ with $dg_i = z_i$ and computed $H_{j-1}( B^{i-1}(D) / B^{i-2}(D) )$, taking note of how the rank was changed.  Finally we accounted for this pesky $\otimes S^{\otimes n-1}$.  

\bigskip

\noindent {\bf Examples.}  Consider the arc presentations of the trefoil $D_T$ and figure-eight knot $D_8$ shown in Figure~\ref{figure:treffig8}. \begin{figure}[htbp]
\begin{center}
\input{treffig8-2.eepic}
\vspace*{-2mm}\caption{Arc presentations $D_T, D_8$ of the trefoil and figure-8 knot}
\label{figure:treffig8}
\end{center}
\vspace*{-5mm}\end{figure}
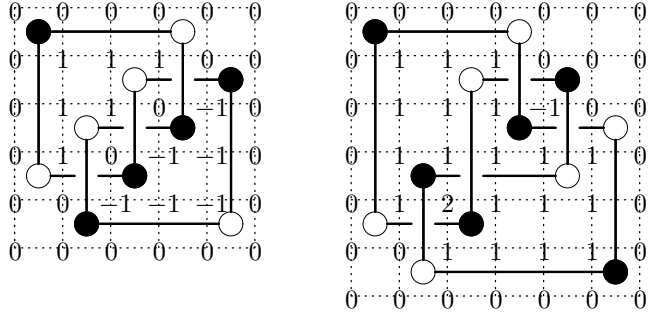 It is easy to see that, among the 120 permutations of $\{1,2,3,4,5\}$, only $51234 \in C_2(D_T)$ and $15234,41234,51243,51324,52134 \in C_1(D_T)$ have non-negative Alexander gradings in $C(D_T)$.  The first of these has Alexander grading 1 while the others have Alexander grading $0$ so there can be no boundary maps between these generators in the complex of successive quotients and we immediately see that $H_j(C^i(D_T) / C^{i-1}(D_T))=C_j^i(D_T) / C_j^{i-1}(D_T)$ for $i \geq 0$, which, together with the symmetry mentioned above yields the correct $\widehat{HFK}$ of the trefoil.  Also notice that $$d(51234) = 15234+41234+51243+51324+52134 \in C(D_T) $$ so it follows that $$d^1:\widehat{HFK}_2(3_1,1) \to \widehat{HFK}_1(3_1,0)$$ has rank 1 as expected.  The situation is similar for $C(D_8)$ where only the generators  $$23146,312456,321465,321546,324156,326451,421356,621453 \in C_0(D_8)$$ and $321456 \in C_1(D_8)$ have non-negative Alexander grading, the latter having Alexander grading $1$ and the others having grading $0$.

\end{section}

\begin{section}{Results} \label{section:results}

We obtained our grid diagrams from several sources.  Originally we used a Mathematica package written by Dror Bar-Natan available at 

\smallskip

{\tt www.math.toronto.edu/$\sim$drorbn/Misc/MOSComplex/index.html}  

\smallskip

\noindent though unfortunately the arc presentations we obtained were generally too large to use directly. We wrote a program to aid in the reduction of arc presentations, though a lot of human labor was still required.  Subsequently, Peter Ozsv\'ath told us about Marc Culler's {\tt gridlink} program \cite{MC}.  We then obtained better diagrams by using his ``simplify" feature which tries making a large number of random grid moves, looking for destabilizations, to reduce the size of a given grid diagram.  This yielded workable diagrams for 12-crossing non-alternating knots as well as much better diagrams for 11-crossing non-alternating knots.  Lenny Ng then got ahold of our list of diagrams for 11-crossing non-alternating knots and checked that all but a half dozen were minimal (by using Ian Nutt's table of knots with arc-index $\leq 10$); he then managed to produce minimal diagrams for the few remaining knots.  Thus we managed to produce a perfect list of grid diagrams for 11-crossing non-alternating knots (available online at our website mentioned at the end of this article).  All such knots have arc-index $\leq 11$.  Similarly we obtained grid diagrams with arc-index at most 13 for all 12-crossing non-alternating knots (Lenny Ng helped us produce an arc-index 12 diagram of $12n_{453}$, which {\tt gridlink} wouldn't get below arc-index 14).  In many cases, we have taken the mirrors of diagrams obtained from {\tt gridlink} (in turn, these are sometimes mirrors of the diagrams from Bar-Natan's program) in order to arrange that the knot Floer homology in Maslov grading 0 is supported in non-negative Alexander grading (this facilitates computation of $\tau$). 

Using the techniques described in Section~\ref{section:methodology} we were able to compute the knot Floer homology for all knots with at most 12 crossings and the $\tau$ invariant for all knots with at most 11 crossings.  Indeed, for knots with at most 11 crossings we actually computed the $E^2$ term of the spectral sequence associated to the knot filtration $\widehat{CFK}$.  Our computations are presented in the tables below.  We do not write down the $E^2$ term if the knot Floer homology is supported on one diagonal.

\newpage

\begin{tiny}

$$
\begin{array}{lll}
{\rm Knot\,} K & \tau & \sum_{i,j} \dim_{\ZZ_2} \widehat{HFK}_j(S^3,K,i) q^jt^i = \sum_{i,j} \dim_{\ZZ_2} E^1_{p,j-p}(\widehat{CFK}(S^3,K)) q^jt^i \\
& & \sum_{i,j} \dim_{\ZZ_2} E^2_{p,j-p}(\widehat{CFK}(S^3,K)) q^jt^i \\
\\
8_{19} & -3 & t^{-3}+qt^{-2}+q^2+q^5t^2+q^6t^3 \\
8_{20} & 0 & q^{-2}t^{-2}+2q^{-1}t^{-1}+3+2qt+q^2t^2 \\
8_{21} & -1 & q^{-1}t^{-2}+4t^{-1}+5q+4q^2t+q^3t^2 \\
9_{42} & 0 & q^{-1}t^{-2}+2t^{-1}+1+2q+2q^2t+q^3t^2 \\
 & & 1\\
\end{array}$$ $$\begin{array}{lll}
11n_{1} & -1 & q^{-1}t^{-2}+7t^{-1}+11q+7q^{2}t+q^{3}t^{2} \\
11n_{2} & -2 & 2q^{-1}t^{-3}+8t^{-2}+12qt^{-1}+13q^{2}+12q^{3}t+8q^{4}t^{2}+2q^{5}t^{3} \\
11n_{3} & -1 & 3q^{-1}t^{-2}+11t^{-1}+15q+11q^{2}t+3q^{3}t^{2} \\
11n_{4} & 0 & q^{-3}t^{-3}+5q^{-2}t^{-2}+11q^{-1}t^{-1}+15+11qt+5q^{2}t^{2}+q^{3}t^{3} \\
11n_{5} & -1 & q^{-2}t^{-3}+7q^{-1}t^{-2}+17t^{-1}+21q+17q^{2}t+7q^{3}t^{2}+q^{4}t^{3} \\
11n_{6} & 0 & q^{-3}t^{-3}+q^{-3}t^{-2}+4q^{-2}t^{-2}+2q^{-2}t^{-1}+5q^{-1}t^{-1}+2q^{-1}+5+2t+5qt+qt^{2}+4q^{2}t^{2}+q^{3}t^{3} \\
 & & 1 \\
11n_{7} & -1 & q^{-2}t^{-3}+6q^{-1}t^{-2}+16t^{-1}+21q+16q^{2}t+6q^{3}t^{2}+q^{4}t^{3} \\
\end{array}$$ $$\begin{array}{lll}
11n_{61} & 2 & q^{-6}t^{-4}+3q^{-5}t^{-3}+4q^{-4}t^{-2}+2q^{-3}t^{-1}+q^{-2}t^{-1}+q^{-1}+2q^{-1}t+t+4t^{2}+3qt^{3}+q^{2}t^{4} \\
 & & q^{-2}t^{-1}+q^{-1}t+t^2 \\
11n_{62} & 0 & 2q^{-2}t^{-2}+8q^{-1}t^{-1}+13+8qt+2q^{2}t^{2} \\
11n_{63} & -1 & 2q^{-1}t^{-2}+10t^{-1}+15q+10q^{2}t+2q^{3}t^{2} \\
11n_{64} & -2 & 2q^{-1}t^{-3}+5t^{-2}+3qt^{-1}+q^{2}+3q^{3}t+5q^{4}t^{2}+2q^{5}t^{3} \\
11n_{65} & 0 & 3q^{-2}t^{-2}+8q^{-1}t^{-1}+11+8qt+3q^{2}t^{2} \\
11n_{66} & -1 & q^{-2}t^{-3}+7q^{-1}t^{-2}+18t^{-1}+23q+18q^{2}t+7q^{3}t^{2}+q^{4}t^{3} \\
11n_{67} & 0 & q^{-3}t^{-2}+q^{-2}t^{-2}+2q^{-2}t^{-1}+4q^{-1}t^{-1}+2q^{-1}+7+2t+4qt+qt^{2}+q^{2}t^{2} \\
\end{array}$$ $$\begin{array}{lll}
11n_{121} & -2 & q^{-1}t^{-3}+6t^{-2}+10qt^{-1}+11q^{2}+10q^{3}t+6q^{4}t^{2}+q^{5}t^{3} \\
11n_{122} & -1 & 2q^{-1}t^{-2}+7t^{-1}+9q+7q^{2}t+2q^{3}t^{2} \\
11n_{123} & 0 & 3q^{-2}t^{-2}+14q^{-1}t^{-1}+23+14qt+3q^{2}t^{2} \\
11n_{124} & -1 & q^{-2}t^{-3}+6q^{-1}t^{-2}+14t^{-1}+17q+14q^{2}t+6q^{3}t^{2}+q^{4}t^{3} \\
11n_{125} & -1 & q^{-2}t^{-3}+6q^{-1}t^{-2}+15t^{-1}+19q+15q^{2}t+6q^{3}t^{2}+q^{4}t^{3} \\
11n_{126} & 3 & 3q^{-6}t^{-3}+6q^{-5}t^{-2}+4q^{-4}t^{-1}+2q^{-3}+q^{-2}+4q^{-2}t+6q^{-1}t^{2}+3t^{3} \\
 & & q^{-2}+q^{-1}t^2+t^3 \\
11n_{127} & -1 & q^{-2}t^{-3}+5q^{-1}t^{-2}+13t^{-1}+17q+13q^{2}t+5q^{3}t^{2}+q^{4}t^{3} \\
\end{array}$$ 

\newpage

$$\begin{array}{ll} 12n_{1} & q^{-3}t^{-3}+5q^{-2}t^{-2}+13q^{-1}t^{-1}+19+13qt+5q^{2}t^{2}+q^{3}t^{3} \\
12n_{2} & q^{-3}t^{-4}+5q^{-2}t^{-3}+11q^{-1}t^{-2}+16t^{-1}+17q+16q^{2}t+11q^{3}t^{2}+5q^{4}t^{3}+q^{5}t^{4} \\
12n_{3} & q^{-3}t^{-3}+7q^{-2}t^{-2}+17q^{-1}t^{-1}+23+17qt+7q^{2}t^{2}+q^{3}t^{3} \\
12n_{4} & q^{-3}t^{-3}+7q^{-2}t^{-2}+19q^{-1}t^{-1}+27+19qt+7q^{2}t^{2}+q^{3}t^{3} \\
12n_{5} & q^{-4}t^{-3}+7q^{-3}t^{-2}+21q^{-2}t^{-1}+29q^{-1}+21t+7qt^{2}+q^{2}t^{3} \\
12n_{6} & q^{-2}t^{-4}+5q^{-1}t^{-3}+9t^{-2}+6qt^{-1}+3q^{2}+6q^{3}t+9q^{4}t^{2}+5q^{5}t^{3}+q^{6}t^{4} \\
12n_{7} & q^{-4}t^{-3}+5q^{-3}t^{-2}+9q^{-2}t^{-1}+10q^{-1}+1+9t+5qt^{2}+q^{2}t^{3} \\
12n_{8} & q^{-2}t^{-4}+5q^{-1}t^{-3}+9t^{-2}+10qt^{-1}+11q^{2}+10q^{3}t+9q^{4}t^{2}+5q^{5}t^{3}+q^{6}t^{4} \\
\end{array}$$ $$\begin{array}{ll}
12n_{66} & 2q^{-3}t^{-3}+9q^{-2}t^{-2}+18q^{-1}t^{-1}+23+18qt+9q^{2}t^{2}+2q^{3}t^{3} \\
12n_{67} & q^{-7}t^{-4}+5q^{-6}t^{-3}+8q^{-5}t^{-2}+8q^{-4}t^{-1}+8q^{-3}+q^{-2}+8q^{-2}t+8q^{-1}t^{2}+5t^{3}+qt^{4} \\
12n_{68} & q^{-7}t^{-4}+5q^{-6}t^{-3}+7q^{-5}t^{-2}+4q^{-4}t^{-1}+2q^{-3}+q^{-2}+4q^{-2}t+7q^{-1}t^{2}+5t^{3}+qt^{4} \\
12n_{69} & q^{-3}t^{-3}+7q^{-2}t^{-2}+17q^{-1}t^{-1}+23+17qt+7q^{2}t^{2}+q^{3}t^{3} \\
12n_{70} & q^{-4}t^{-3}+7q^{-3}t^{-2}+19q^{-2}t^{-1}+25q^{-1}+19t+7qt^{2}+q^{2}t^{3} \\
12n_{71} & q^{-5}t^{-4}+5q^{-4}t^{-3}+9q^{-3}t^{-2}+10q^{-2}t^{-1}+9q^{-1}+10t+9qt^{2}+5q^{2}t^{3}+q^{3}t^{4} \\
12n_{72} & q^{-5}t^{-3}+9q^{-4}t^{-2}+23q^{-3}t^{-1}+31q^{-2}+23q^{-1}t+9t^{2}+qt^{3} \\
12n_{73} & q^{-4}t^{-3}+5q^{-3}t^{-2}+13q^{-2}t^{-1}+17q^{-1}+13t+5qt^{2}+q^{2}t^{3} \\
\end{array}$$ $$\begin{array}{ll}
12n_{131} & 2q^{-4}t^{-3}+10q^{-3}t^{-2}+24q^{-2}t^{-1}+31q^{-1}+24t+10qt^{2}+2q^{2}t^{3} \\
12n_{132} & q^{-4}t^{-3}+q^{-3}t^{-3}+2q^{-3}t^{-2}+6q^{-2}t^{-2}+q^{-2}t^{-1}+13q^{-1}t^{-1}+17+t+13qt+2qt^{2}+6q^{2}t^{2}+q^{2}t^{3}+q^{3}t^{3} \\
12n_{133} & 4t^{-3}+14qt^{-2}+26q^{2}t^{-1}+31q^{3}+26q^{4}t+14q^{5}t^{2}+4q^{6}t^{3} \\
12n_{134} & q^{-7}t^{-4}+4q^{-6}t^{-3}+6q^{-5}t^{-2}+7q^{-4}t^{-1}+8q^{-3}+q^{-2}+7q^{-2}t+6q^{-1}t^{2}+4t^{3}+qt^{4} \\
12n_{135} & q^{-2}t^{-4}+4q^{-1}t^{-3}+10t^{-2}+15qt^{-1}+17q^{2}+15q^{3}t+10q^{4}t^{2}+4q^{5}t^{3}+q^{6}t^{4} \\
12n_{136} & q^{-8}t^{-4}+q^{-7}t^{-3}+q^{-6}t^{-3}+6q^{-5}t^{-2}+15q^{-4}t^{-1}+19q^{-3}+15q^{-2}t+6q^{-1}t^{2}+q^{-1}t^{3}+t^{3}+t^{4} \\
12n_{137} & q^{-4}t^{-4}+4q^{-3}t^{-3}+12q^{-2}t^{-2}+23q^{-1}t^{-1}+29+23qt+12q^{2}t^{2}+4q^{3}t^{3}+q^{4}t^{4} \\
12n_{138} & q^{-7}t^{-4}+4q^{-6}t^{-3}+5q^{-5}t^{-2}+q^{-4}t^{-2}+3q^{-4}t^{-1}+4q^{-3}t^{-1}+2q^{-3}+7q^{-2}+3q^{-2}t+4q^{-1}t+5q^{-1}t^{2}+t^{2}+4t^{3}+qt^{4} \\
\end{array}$$ $$\begin{array}{ll}
12n_{197} & 2q^{-2}t^{-3}+10q^{-1}t^{-2}+25t^{-1}+33q+25q^{2}t+10q^{3}t^{2}+2q^{4}t^{3} \\
12n_{198} & q^{-2}t^{-3}+4q^{-1}t^{-2}+q^{-1}t^{-1}+7t^{-1}+3+8q+qt+7q^{2}t+4q^{3}t^{2}+q^{4}t^{3} \\
12n_{199} & q^{-4}t^{-3}+q^{-4}t^{-2}+2q^{-3}t^{-2}+q^{-3}t^{-1}+3q^{-2}t^{-1}+q^{-2}+4q^{-1}+q^{-1}t+3t+t^{2}+2qt^{2}+q^{2}t^{3} \\
12n_{200} & 2q^{-2}t^{-2}+4q^{-1}t^{-1}+2t^{-1}+4+3q+4qt+2q^{2}t+2q^{2}t^{2} \\
12n_{201} & 3q^{-5}t^{-3}+13q^{-4}t^{-2}+23q^{-3}t^{-1}+27q^{-2}+23q^{-1}t+13t^{2}+3qt^{3} \\
12n_{202} & 3q^{-2}t^{-3}+13q^{-1}t^{-2}+31t^{-1}+41q+31q^{2}t+13q^{3}t^{2}+3q^{4}t^{3} \\
12n_{203} & 3q^{-6}t^{-3}+6q^{-5}t^{-2}+q^{-4}t^{-2}+5q^{-4}t^{-1}+4q^{-3}t^{-1}+3q^{-3}+6q^{-2}+5q^{-2}t+4q^{-1}t+6q^{-1}t^{2}+t^{2}+3t^{3} \\
12n_{204} & q^{-5}t^{-3}+9q^{-4}t^{-2}+22q^{-3}t^{-1}+29q^{-2}+22q^{-1}t+9t^{2}+qt^{3} \\
\end{array}$$ $$\begin{array}{ll}
12n_{266} & 2q^{-3}t^{-3}+10q^{-2}t^{-2}+24q^{-1}t^{-1}+33+24qt+10q^{2}t^{2}+2q^{3}t^{3} \\
12n_{267} & q^{-5}t^{-3}+q^{-4}t^{-3}+2q^{-4}t^{-2}+6q^{-3}t^{-2}+q^{-3}t^{-1}+17q^{-2}t^{-1}+23q^{-1}+q^{-1}t+17t+2t^{2}+6qt^{2}+qt^{3}+q^{2}t^{3} \\
12n_{268} & q^{-3}t^{-2}+q^{-2}t^{-2}+4q^{-2}t^{-1}+6q^{-1}t^{-1}+6q^{-1}+11+4t+6qt+qt^{2}+q^{2}t^{2} \\
12n_{269} & q^{-3}t^{-3}+9q^{-2}t^{-2}+28q^{-1}t^{-1}+41+28qt+9q^{2}t^{2}+q^{3}t^{3} \\
12n_{270} & 3q^{-3}t^{-2}+16q^{-2}t^{-1}+25q^{-1}+16t+3qt^{2} \\
12n_{271} & q^{-3}t^{-3}+5q^{-2}t^{-2}+12q^{-1}t^{-1}+17+12qt+5q^{2}t^{2}+q^{3}t^{3} \\
12n_{272} & q^{-4}t^{-3}+7q^{-3}t^{-2}+20q^{-2}t^{-1}+27q^{-1}+20t+7qt^{2}+q^{2}t^{3} \\
12n_{273} & q^{-5}t^{-3}+3q^{-4}t^{-2}+3q^{-3}t^{-1}+3q^{-2}t^{-1}+2q^{-2}+5q^{-1}+3q^{-1}t+3t+3t^{2}+qt^{3} \\
\end{array}$$ $$\begin{array}{ll}
12n_{331} & 2q^{-5}t^{-3}+7q^{-4}t^{-2}+9q^{-3}t^{-1}+9q^{-2}+9q^{-1}t+7t^{2}+2qt^{3} \\
12n_{332} & 2q^{-4}t^{-2}+4q^{-3}t^{-1}+2q^{-2}t^{-1}+5q^{-2}+4q^{-1}+4q^{-1}t+2t+2t^{2} \\
12n_{333} & 3q^{-2}t^{-2}+11q^{-1}t^{-1}+17+11qt+3q^{2}t^{2} \\
12n_{334} & 4q^{-1}t^{-2}+16t^{-1}+23q+16q^{2}t+4q^{3}t^{2} \\
12n_{335} & 2q^{-4}t^{-3}+10q^{-3}t^{-2}+22q^{-2}t^{-1}+27q^{-1}+22t+10qt^{2}+2q^{2}t^{3} \\
12n_{336} & q^{-5}t^{-3}+3q^{-4}t^{-2}+q^{-3}t^{-2}+5q^{-3}t^{-1}+5q^{-2}t^{-1}+6q^{-2}+7q^{-1}+5q^{-1}t+5t+3t^{2}+qt^{2}+qt^{3} \\
12n_{337} & q^{-2}t^{-3}+9q^{-1}t^{-2}+25t^{-1}+33q+25q^{2}t+9q^{3}t^{2}+q^{4}t^{3} \\
12n_{338} & 2q^{-8}t^{-4}+4q^{-7}t^{-3}+4q^{-6}t^{-2}+2q^{-5}t^{-1}+q^{-4}+2q^{-3}t+4q^{-2}t^{2}+4q^{-1}t^{3}+2t^{4} \\
\end{array}$$ $$\begin{array}{ll}
12n_{399} & q^{-3}t^{-3}+7q^{-2}t^{-2}+19q^{-1}t^{-1}+27+19qt+7q^{2}t^{2}+q^{3}t^{3} \\
12n_{400} & q^{-4}t^{-3}+8q^{-3}t^{-2}+23q^{-2}t^{-1}+31q^{-1}+23t+8qt^{2}+q^{2}t^{3} \\
12n_{401} & q^{-2}t^{-3}+7q^{-1}t^{-2}+21t^{-1}+29q+21q^{2}t+7q^{3}t^{2}+q^{4}t^{3} \\
12n_{402} & t^{-3}+2t^{-2}+qt^{-2}+4qt^{-1}+q^{2}t^{-1}+4q^{2}+q^{3}+4q^{3}t+q^{4}t+2q^{4}t^{2}+q^{5}t^{2}+q^{6}t^{3} \\
12n_{403} & q^{-1}t^{-3}+3t^{-2}+2t^{-1}+3qt^{-1}+3q+2q^{2}+2q^{2}t+3q^{3}t+3q^{4}t^{2}+q^{5}t^{3} \\
12n_{404} & 2t^{-2}+4t^{-1}+3qt^{-1}+8q+3q^{2}+4q^{2}t+3q^{3}t+2q^{4}t^{2} \\
12n_{405} & q^{-6}t^{-4}+6q^{-5}t^{-3}+13q^{-4}t^{-2}+12q^{-3}t^{-1}+9q^{-2}+12q^{-1}t+13t^{2}+6qt^{3}+q^{2}t^{4} \\
12n_{406} & 5q^{-6}t^{-3}+13q^{-5}t^{-2}+16q^{-4}t^{-1}+15q^{-3}+16q^{-2}t+13q^{-1}t^{2}+5t^{3} \\
\end{array}$$ $$\begin{array}{ll}
12n_{468} & q^{-4}t^{-4}+3q^{-3}t^{-3}+5q^{-2}t^{-2}+6q^{-1}t^{-1}+7+6qt+5q^{2}t^{2}+3q^{3}t^{3}+q^{4}t^{4} \\
12n_{469} & 2q^{-2}t^{-3}+6q^{-1}t^{-2}+12t^{-1}+15q+12q^{2}t+6q^{3}t^{2}+2q^{4}t^{3} \\
12n_{470} & 2q^{-3}t^{-3}+6q^{-2}t^{-2}+8q^{-1}t^{-1}+9+8qt+6q^{2}t^{2}+2q^{3}t^{3} \\
12n_{471} & 2q^{-4}t^{-3}+8q^{-3}t^{-2}+16q^{-2}t^{-1}+19q^{-1}+16t+8qt^{2}+2q^{2}t^{3} \\
12n_{472} & q^{-10}t^{-5}+q^{-9}t^{-4}+2q^{-6}t^{-2}+4q^{-5}t^{-1}+5q^{-4}+4q^{-3}t+2q^{-2}t^{2}+q^{-1}t^{4}+t^{5} \\
12n_{473} & 2q^{-8}t^{-4}+3q^{-7}t^{-3}+q^{-6}t^{-2}+3q^{-4}t^{-1}+5q^{-3}+3q^{-2}t+q^{-2}t^{2}+3q^{-1}t^{3}+2t^{4} \\
12n_{474} & 2q^{-8}t^{-4}+3q^{-7}t^{-3}+2q^{-6}t^{-2}+q^{-5}t^{-1}+q^{-4}+q^{-3}t+2q^{-2}t^{2}+3q^{-1}t^{3}+2t^{4} \\
12n_{475} & q^{-2}t^{-3}+2q^{-1}t^{-2}+t^{-1}+1+q^{2}t+2q^{3}t^{2}+q^{4}t^{3} \\
\end{array}$$ $$\begin{array}{ll}
12n_{526} & q^{-2}t^{-4}+4q^{-1}t^{-3}+q^{-1}t^{-2}+7t^{-2}+2t^{-1}+5qt^{-1}+2q+3q^{2}+2q^{2}t+5q^{3}t+q^{3}t^{2}+7q^{4}t^{2}+4q^{5}t^{3}+q^{6}t^{4} \\
12n_{527} & 3q^{-1}t^{-3}+12t^{-2}+20qt^{-1}+23q^{2}+20q^{3}t+12q^{4}t^{2}+3q^{5}t^{3} \\
12n_{528} & q^{-6}t^{-3}+q^{-5}t^{-2}+3q^{-4}t^{-2}+10q^{-3}t^{-1}+15q^{-2}+10q^{-1}t+q^{-1}t^{2}+3t^{2}+t^{3} \\
12n_{529} & q^{-2}t^{-3}+10q^{-1}t^{-2}+28t^{-1}+37q+28q^{2}t+10q^{3}t^{2}+q^{4}t^{3} \\
12n_{530} & 2q^{-2}t^{-3}+10q^{-1}t^{-2}+19t^{-1}+21q+19q^{2}t+10q^{3}t^{2}+2q^{4}t^{3} \\
12n_{531} & q^{-4}t^{-4}+4q^{-3}t^{-3}+11q^{-2}t^{-2}+20q^{-1}t^{-1}+25+20qt+11q^{2}t^{2}+4q^{3}t^{3}+q^{4}t^{4} \\
12n_{532} & 2q^{-3}t^{-3}+11q^{-2}t^{-2}+29q^{-1}t^{-1}+41+29qt+11q^{2}t^{2}+2q^{3}t^{3} \\
12n_{533} & q^{-5}t^{-4}+4q^{-4}t^{-3}+8q^{-3}t^{-2}+12q^{-2}t^{-1}+13q^{-1}+12t+8qt^{2}+4q^{2}t^{3}+q^{3}t^{4} \\
\end{array}$$ $$\begin{array}{ll}
12n_{592} & q^{-6}t^{-4}+4q^{-5}t^{-3}+8q^{-4}t^{-2}+8q^{-3}t^{-1}+7q^{-2}+8q^{-1}t+8t^{2}+4qt^{3}+q^{2}t^{4} \\
12n_{593} & 2q^{-5}t^{-3}+9q^{-4}t^{-2}+13q^{-3}t^{-1}+13q^{-2}+13q^{-1}t+9t^{2}+2qt^{3} \\
12n_{594} & 2t^{-3}+2t^{-2}+3qt^{-2}+6qt^{-1}+q^{2}t^{-1}+9q^{2}+6q^{3}t+q^{4}t+2q^{4}t^{2}+3q^{5}t^{2}+2q^{6}t^{3} \\
12n_{595} & 2q^{-4}t^{-3}+9q^{-3}t^{-2}+21q^{-2}t^{-1}+27q^{-1}+21t+9qt^{2}+2q^{2}t^{3} \\
12n_{596} & 3q^{-2}t^{-3}+11q^{-1}t^{-2}+21t^{-1}+25q+21q^{2}t+11q^{3}t^{2}+3q^{4}t^{3} \\
12n_{597} & 2q^{-2}t^{-3}+8q^{-1}t^{-2}+17t^{-1}+21q+17q^{2}t+8q^{3}t^{2}+2q^{4}t^{3} \\
12n_{598} & q^{-2}t^{-3}+8q^{-1}t^{-2}+24t^{-1}+33q+24q^{2}t+8q^{3}t^{2}+q^{4}t^{3} \\
12n_{599} & q^{-3}t^{-4}+4q^{-2}t^{-3}+7q^{-1}t^{-2}+10t^{-1}+11q+10q^{2}t+7q^{3}t^{2}+4q^{4}t^{3}+q^{5}t^{4} \\
\end{array}$$ $$\begin{array}{ll}
12n_{657} & q^{-4}t^{-4}+4q^{-3}t^{-3}+9q^{-2}t^{-2}+16q^{-1}t^{-1}+21+16qt+9q^{2}t^{2}+4q^{3}t^{3}+q^{4}t^{4} \\
12n_{658} & q^{-6}t^{-4}+5q^{-5}t^{-3}+12q^{-4}t^{-2}+17q^{-3}t^{-1}+19q^{-2}+17q^{-1}t+12t^{2}+5qt^{3}+q^{2}t^{4} \\
12n_{659} & 2q^{-1}t^{-3}+10t^{-2}+16qt^{-1}+17q^{2}+16q^{3}t+10q^{4}t^{2}+2q^{5}t^{3} \\
12n_{660} & t^{-3}+4t^{-2}+qt^{-2}+12qt^{-1}+17q^{2}+12q^{3}t+4q^{4}t^{2}+q^{5}t^{2}+q^{6}t^{3} \\
12n_{661} & 2q^{-1}t^{-3}+6t^{-2}+5qt^{-1}+3q^{2}+5q^{3}t+6q^{4}t^{2}+2q^{5}t^{3} \\
12n_{662} & 2q^{-2}t^{-3}+8q^{-1}t^{-2}+13t^{-1}+13q+13q^{2}t+8q^{3}t^{2}+2q^{4}t^{3} \\
12n_{663} & q^{-3}t^{-4}+6q^{-2}t^{-3}+15q^{-1}t^{-2}+23t^{-1}+25q+23q^{2}t+15q^{3}t^{2}+6q^{4}t^{3}+q^{5}t^{4} \\
12n_{664} & q^{-4}t^{-4}+7q^{-3}t^{-3}+20q^{-2}t^{-2}+34q^{-1}t^{-1}+41+34qt+20q^{2}t^{2}+7q^{3}t^{3}+q^{4}t^{4} \\
\end{array}$$ $$\begin{array}{ll}
12n_{726} & 4q^{-1}t^{-2}+16t^{-1}+23q+16q^{2}t+4q^{3}t^{2} \\
12n_{727} & 2q^{-2}t^{-3}+10q^{-1}t^{-2}+21t^{-1}+25q+21q^{2}t+10q^{3}t^{2}+2q^{4}t^{3} \\
12n_{728} & q^{-3}t^{-4}+5q^{-2}t^{-3}+13q^{-1}t^{-2}+22t^{-1}+25q+22q^{2}t+13q^{3}t^{2}+5q^{4}t^{3}+q^{5}t^{4} \\
12n_{729} & 2q^{-1}t^{-3}+9t^{-2}+16qt^{-1}+19q^{2}+16q^{3}t+9q^{4}t^{2}+2q^{5}t^{3} \\
12n_{730} & q^{-2}t^{-3}+4q^{-1}t^{-2}+11t^{-1}+15q+11q^{2}t+4q^{3}t^{2}+q^{4}t^{3} \\
12n_{731} & 2q^{-3}t^{-3}+11q^{-2}t^{-2}+28q^{-1}t^{-1}+39+28qt+11q^{2}t^{2}+2q^{3}t^{3} \\
12n_{732} & q^{-2}t^{-3}+9q^{-1}t^{-2}+28t^{-1}+39q+28q^{2}t+9q^{3}t^{2}+q^{4}t^{3} \\
12n_{733} & 2q^{-4}t^{-3}+12q^{-3}t^{-2}+31q^{-2}t^{-1}+41q^{-1}+31t+12qt^{2}+2q^{2}t^{3} \\
\end{array}$$ $$\begin{array}{ll}
12n_{791} & q^{-2}t^{-4}+5q^{-1}t^{-3}+13t^{-2}+18qt^{-1}+19q^{2}+18q^{3}t+13q^{4}t^{2}+5q^{5}t^{3}+q^{6}t^{4} \\
12n_{792} & q^{-3}t^{-4}+4q^{-2}t^{-3}+10q^{-1}t^{-2}+16t^{-1}+17q+16q^{2}t+10q^{3}t^{2}+4q^{4}t^{3}+q^{5}t^{4} \\
12n_{793} & 3q^{-1}t^{-3}+9t^{-2}+12qt^{-1}+13q^{2}+12q^{3}t+9q^{4}t^{2}+3q^{5}t^{3} \\
12n_{794} & 2q^{-3}t^{-3}+9q^{-2}t^{-2}+23q^{-1}t^{-1}+33+23qt+9q^{2}t^{2}+2q^{3}t^{3} \\
12n_{795} & 6q^{-1}t^{-2}+20t^{-1}+27q+20q^{2}t+6q^{3}t^{2} \\
12n_{796} & 2q^{-1}t^{-3}+11t^{-2}+22qt^{-1}+27q^{2}+22q^{3}t+11q^{4}t^{2}+2q^{5}t^{3} \\
12n_{797} & 3q^{-1}t^{-2}+15t^{-1}+23q+15q^{2}t+3q^{3}t^{2} \\
12n_{798} & q^{-2}t^{-3}+7q^{-1}t^{-2}+22t^{-1}+31q+22q^{2}t+7q^{3}t^{2}+q^{4}t^{3} \\
\end{array}$$ $$\begin{array}{ll}
12n_{860} & 3q^{-4}t^{-3}+15q^{-3}t^{-2}+33q^{-2}t^{-1}+41q^{-1}+33t+15qt^{2}+3q^{2}t^{3} \\
12n_{861} & 2q^{-2}t^{-3}+9q^{-1}t^{-2}+22t^{-1}+29q+22q^{2}t+9q^{3}t^{2}+2q^{4}t^{3} \\
12n_{862} & 3q^{-2}t^{-3}+13q^{-1}t^{-2}+29t^{-1}+37q+29q^{2}t+13q^{3}t^{2}+3q^{4}t^{3} \\
12n_{863} & 3q^{-1}t^{-3}+16t^{-2}+33qt^{-1}+41q^{2}+33q^{3}t+16q^{4}t^{2}+3q^{5}t^{3} \\
12n_{864} & q^{-5}t^{-4}+5q^{-4}t^{-3}+15q^{-3}t^{-2}+28q^{-2}t^{-1}+33q^{-1}+28t+15qt^{2}+5q^{2}t^{3}+q^{3}t^{4} \\
12n_{865} & q^{-3}t^{-4}+4q^{-2}t^{-3}+11q^{-1}t^{-2}+21t^{-1}+25q+21q^{2}t+11q^{3}t^{2}+4q^{4}t^{3}+q^{5}t^{4} \\
12n_{866} & q^{-4}t^{-4}+5q^{-3}t^{-3}+15q^{-2}t^{-2}+30q^{-1}t^{-1}+39+30qt+15q^{2}t^{2}+5q^{3}t^{3}+q^{4}t^{4} \\
12n_{867} & 2q^{-1}t^{-3}+14t^{-2}+31qt^{-1}+39q^{2}+31q^{3}t+14q^{4}t^{2}+2q^{5}t^{3} \\
\end{array}$$ 

\end{tiny}

\end{section}

\newpage

\begin{section}{Discussion} \label{section:discussion}

The knots $11n_{42}$ and $11n_{34}$ are often discussed in the literature. These are the Kinoshita-Terasaka knot $KT_{2,1}$ and its Conway mutant $C_{2,1}$, respectively. A minimal\footnote{Ian Nutt listed all knots with arc-index at most 10 in his Ph.D. thesis \cite{N}.  $KT_{2,1}$ is not among them.} arc presentation for $KT_{2,1}$ appears in Figure~\ref{figure:kino} from which we computed its knot Floer homology polynomial: \begin{small} $$(q^{-2} + q^{-1})t^{-2} + 4(q^{-1} + 1)t^{-1} + 7 + 6q + 4(q + q^{2})t + (q^{2} + q^{3})t^{2}$$ \end{small}  Likewise, the knot Floer homology polynomial of $C_{2,1}$ is: \begin{small} $$(q^{-3} + q^{-2})t^{-3} + 3(q^{-2} + q^{-1})t^{-2} + 3(q^{-1} + 1)t^{-1} + 3 + 2q + 3(q + q^{2})t + 3(q^{2} + q^{3})t^{2} + (q^{3} + q^{4})t^{3}$$ \end{small}  Recall that these knots are of special interest because their Alexander polynomials are trivial. They can be distinguished by their Seifert genera, however: $KT_{2,1}$ has genus 2, whereas $C_{2,1}$ has genus 3 (c.f. \cite{OSz1}, Theorem~1.2).  The homology in extremal Alexander gradings was computed for these knots in \cite{OSz2}, Theorem~1.1.

Regarding the $\tau$ invariant, we found that for knots through 11 crossings, the $E^1$ term of the spectral sequence associated to $\widehat{CFK}$ (i.e. the knot Floer homology) was supported along (at most) two diagonals, one with Euler characteristic $0$ and the other with Euler characteristic $1$.  The simplest possible behaviour for the $d^1$ differential would be to have $E^2=H_\bullet(E^1,d^1)=0$ along the former diagonal and $\dim_{\ZZ_2} E^2 = 1$ supported in Maslov grading $0$.  This is what happens most of the time, but not always.  For example, this is not possible if the Euler characteristic $1$ diagonal is not supported in Maslov grading $0$, as is the case for $11n_{81}$ where the $E^2$ term can be computed purely for ``shape" reasons.  Sometimes the $d^1$ differential is more unpredictable as we see for $10_{154}$, where the $E^1$ and $E^2$ terms and the ranks of $d^1$ and $d^2$ are shown in Figure~\ref{figure:10_154}.  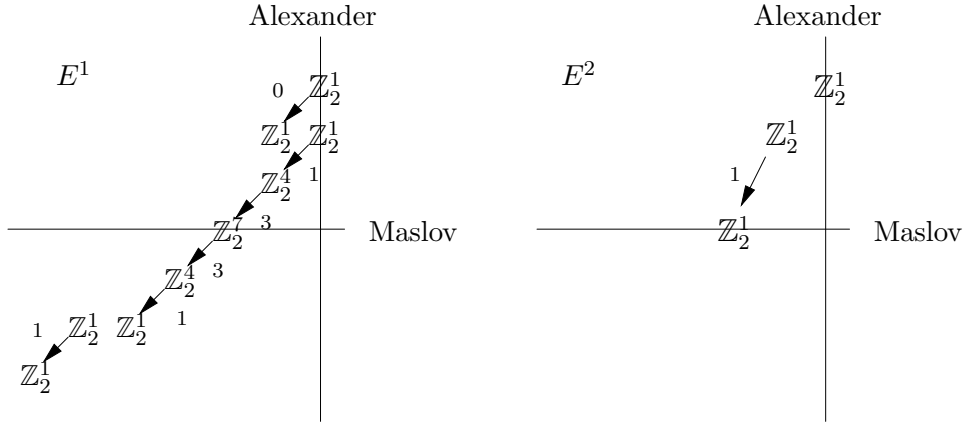
\begin{figure}[htbp]
\begin{center}
\input{10_154.eepic}
\vspace*{-2mm}\caption{The spectral sequence associated to $\widehat{CFK}(10_{154})$}
\label{figure:10_154}
\end{center}
\vspace*{-5mm}\end{figure} We did find that the Rasmussen $s$ invariant is equal to $2 \tau$ for all knots whose $\tau$ invariant we computed, despite the fact that they disagree in general (see \cite{HO} for discussion).

The fact that $s$ and $\tau$ disagree can now be checked by computer, using the examples in \cite{HO}.  We calculated $\widehat{HFK}$ using {\tt gridlink} and our program for the $t$-twisted Whitehead double of the right-handed trefoil for $t=0,\dots,6$.  The results (below) agree with the cases mentioned in \cite{HO} and can also be deduced from the general results of M. Hedden \cite{MH}.

$$\begin{array}{ll} t=0 & 2q^{-3}t^{-1}+2q^{-2}t^{-1}+4q^{-2}+3q^{-1}+2q^{-1}t+2t \\ 
t=1 & 2q^{-3}t^{-1}+q^{-2}t^{-1}+4q^{-2}+q^{-1}+2q^{-1}t+t \\ 
t=2 & 2q^{-3}t^{-1}+4q^{-2}+1+2q^{-1}t\\
t=3 & 2q^{-3}t^{-1}+q^{-1}t^{-1}+4q^{-2}+3+2q^{-1}t+qt \\
t=4 & 2q^{-3}t^{-1}+2q^{-1}t^{-1}+4q^{-2}+5+2q^{-1}t+2qt \\
t=5 & 2q^{-3}t^{-1}+3q^{-1}t^{-1}+4q^{-2}+7+2q^{-1}t+3qt \\
t=6 & 2q^{-3}t^{-1}+4q^{-1}t^{-1}+4q^{-2}+9+2q^{-1}t+4qt \\
\end{array}$$

We conclude with a few comments about arc-index.  There are a surprisingly large number of knots with arc-index at most 13 (e.g. all non-alternating knots with at most 12 crossings).  Non-alternating knots tend to have lower arc-index than alternating ones.  For example, $\alpha(10_{124})=8$ but Cromwell \cite{C} showed that all alternating knots $K$ with 10 or fewer crossings have $\alpha(K) = c(K)+2$ (this fact is proved for alternating knots with any number of crossings in \cite{BP}, which is fortunate since $\widehat{HFK}$ is determined by more easily computed invariants (Alexander polynomial and signature) in the case of alternating knots \cite{OSz4}.  Ian Nutt's table of arc-index $9$ knots mentions two knots with 12 crossings.

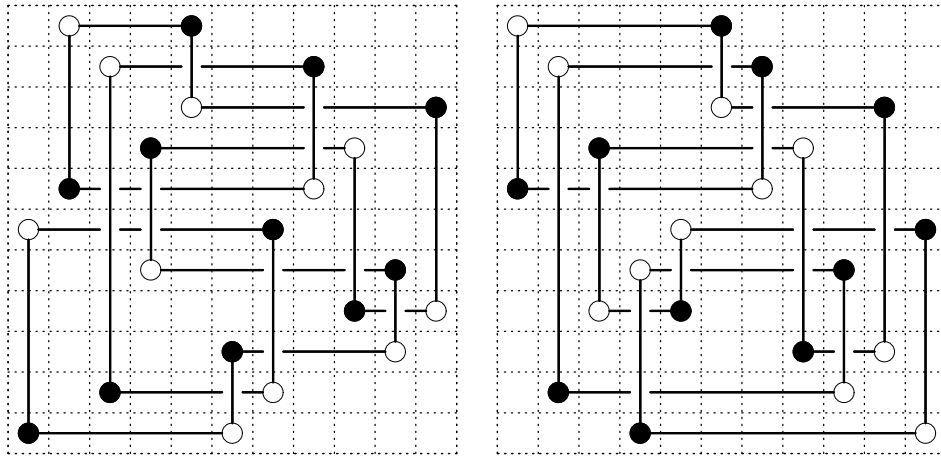
\begin{figure}[htbp]
\begin{center}
\input{kino.eepic}
\vspace*{-2mm}\caption{Arc presentations for $KT_{2,1}$ and $C_{2,1}$}
\label{figure:kino}
\end{center}
\vspace*{-5mm}\end{figure}

Our computer program is available in C++ at the address:
\smallskip

{\tt http://www.math.columbia.edu/$\sim$wgillam/hfk }

\smallskip

We wish to thank Peter Ozsv\'ath, from whom we learned about Heegaard Floer homology, for his encouragement.  We also express our gratitude to Marc Culler for helping us improve our program and for writing {\tt gridlink}, without which we could not have carried out our computations for 12-crossing knots, and to Lenny Ng who drew many minimal diagrams for us by hand. Finally, we thank Josh Greene for his helpful comments on the first draft of this paper.

\end{section}

\end{document}

%% file: rect.eepic
\setlength{\unitlength}{0.00083333in}
\begingroup\makeatletter\ifx\SetFigFont\undefined%
\gdef\SetFigFont#1#2#3#4#5{%
  \reset@font\fontsize{#1}{#2pt}%
  \fontfamily{#3}\fontseries{#4}\fontshape{#5}%
  \selectfont}%
\fi\endgroup%
{\renewcommand{\dashlinestretch}{30}
\begin{picture}(1707,1855)(0,-10)
\path(42,1798)(42,298)
\path(42,298)(1542,298)
\path(1542,298)(1542,1798)(42,1798)
\path(492,298)(492,1798)
\path(1092,298)(1092,1798)
\path(42,1348)(1542,1348)
\path(42,748)(1542,748)
\dottedline{45}(267,298)(267,1798)
\dottedline{45}(792,298)(792,1798)
\dottedline{45}(1317,298)(1317,1798)
\dottedline{45}(42,1573)(1542,1573)
\dottedline{45}(42,1048)(1542,1048)
\dottedline{45}(42,523)(1542,523)
\dottedline{45}(717,1798)(792,1798)
\dottedline{45}(642,1798)(867,1798)
\dottedline{45}(642,1798)(867,1798)
\blacken\path(747.000,1768.000)(867.000,1798.000)(747.000,1828.000)(747.000,1768.000)
\dottedline{45}(642,298)(867,298)
\blacken\path(747.000,268.000)(867.000,298.000)(747.000,328.000)(747.000,268.000)
\dottedline{45}(42,823)(42,973)
\blacken\path(72.000,853.000)(42.000,973.000)(12.000,853.000)(72.000,853.000)
\dottedline{45}(42,973)(42,1198)
\blacken\path(72.000,1078.000)(42.000,1198.000)(12.000,1078.000)(72.000,1078.000)
\dottedline{45}(1542,748)(1542,973)
\blacken\path(1572.000,853.000)(1542.000,973.000)(1512.000,853.000)(1572.000,853.000)
\dottedline{45}(1542,1048)(1542,1198)
\blacken\path(1572.000,1078.000)(1542.000,1198.000)(1512.000,1078.000)(1572.000,1078.000)
\put(717,448){\makebox(0,0)[lb]{\smash{{\SetFigFont{12}{14.4}{\familydefault}{\mddefault}{\updefault}$R_4$}}}}
\put(117,973){\makebox(0,0)[lb]{\smash{{\SetFigFont{12}{14.4}{\familydefault}{\mddefault}{\updefault}$R_3$}}}}
\put(717,1498){\makebox(0,0)[lb]{\smash{{\SetFigFont{12}{14.4}{\familydefault}{\mddefault}{\updefault}$R_4$}}}}
\put(417,73){\makebox(0,0)[lb]{\smash{{\SetFigFont{12}{14.4}{\familydefault}{\mddefault}{\updefault}$\alpha_i$}}}}
\put(1017,73){\makebox(0,0)[lb]{\smash{{\SetFigFont{12}{14.4}{\familydefault}{\mddefault}{\updefault}$\alpha_j$}}}}
\put(1167,973){\makebox(0,0)[lb]{\smash{{\SetFigFont{12}{14.4}{\familydefault}{\mddefault}{\updefault}$R_3$}}}}
\put(1692,673){\makebox(0,0)[lb]{\smash{{\SetFigFont{12}{14.4}{\familydefault}{\mddefault}{\updefault}$\beta_{\sigma(i)}$}}}}
\put(1692,1273){\makebox(0,0)[lb]{\smash{{\SetFigFont{12}{14.4}{\familydefault}{\mddefault}{\updefault}$\beta_{\sigma(j)}$}}}}
\put(717,973){\makebox(0,0)[lb]{\smash{{\SetFigFont{12}{14.4}{\familydefault}{\mddefault}{\updefault}$R_1$}}}}
\put(117,1498){\makebox(0,0)[lb]{\smash{{\SetFigFont{12}{14.4}{\familydefault}{\mddefault}{\updefault}$R_2$}}}}
\put(117,448){\makebox(0,0)[lb]{\smash{{\SetFigFont{12}{14.4}{\familydefault}{\mddefault}{\updefault}$R_2$}}}}
\put(1167,448){\makebox(0,0)[lb]{\smash{{\SetFigFont{12}{14.4}{\familydefault}{\mddefault}{\updefault}$R_2$}}}}
\put(1167,1498){\makebox(0,0)[lb]{\smash{{\SetFigFont{12}{14.4}{\familydefault}{\mddefault}{\updefault}$R_2$}}}}
\end{picture}
}

%% file: chain.eepic
\setlength{\unitlength}{0.00083333in}
\begingroup\makeatletter\ifx\SetFigFont\undefined%
\gdef\SetFigFont#1#2#3#4#5{%
  \reset@font\fontsize{#1}{#2pt}%
  \fontfamily{#3}\fontseries{#4}\fontshape{#5}%
  \selectfont}%
\fi\endgroup%
{\renewcommand{\dashlinestretch}{30}
\begin{picture}(5471,1609)(0,-10)
\put(3165,1273){\blacken\ellipse{96}{96}}
\put(3165,1273){\ellipse{96}{96}}
\put(3615,373){\blacken\ellipse{96}{96}}
\put(3615,373){\ellipse{96}{96}}
\put(4065,1273){\blacken\ellipse{96}{96}}
\put(4065,1273){\ellipse{96}{96}}
\put(4515,373){\blacken\ellipse{96}{96}}
\put(4515,373){\ellipse{96}{96}}
\put(4965,373){\blacken\ellipse{96}{96}}
\put(4965,373){\ellipse{96}{96}}
\put(5415,373){\blacken\ellipse{96}{96}}
\put(5415,373){\ellipse{96}{96}}
\put(90,1273){\blacken\ellipse{96}{96}}
\put(90,1273){\ellipse{96}{96}}
\put(540,373){\blacken\ellipse{96}{96}}
\put(540,373){\ellipse{96}{96}}
\put(990,1273){\blacken\ellipse{96}{96}}
\put(990,1273){\ellipse{96}{96}}
\put(990,373){\blacken\ellipse{96}{96}}
\put(990,373){\ellipse{96}{96}}
\put(1440,373){\blacken\ellipse{96}{96}}
\put(1440,373){\ellipse{96}{96}}
\put(1890,373){\blacken\ellipse{96}{96}}
\put(1890,373){\ellipse{96}{96}}
\put(1890,1273){\blacken\ellipse{96}{96}}
\put(1890,1273){\ellipse{96}{96}}
\put(2340,373){\blacken\ellipse{96}{96}}
\put(2340,373){\ellipse{96}{96}}
\path(2490,973)(2940,973)
\blacken\path(2820.000,943.000)(2940.000,973.000)(2820.000,1003.000)(2820.000,943.000)
\path(3165,1273)(3615,373)
\path(3615,373)(4065,1273)
\dottedline{45}(3315,973)(3465,673)
\blacken\path(3384.502,766.915)(3465.000,673.000)(3438.167,793.748)(3384.502,766.915)
\dottedline{45}(3915,973)(3765,673)
\blacken\path(3791.833,793.748)(3765.000,673.000)(3845.498,766.915)(3791.833,793.748)
\path(4065,1273)(5415,373)
\path(90,1273)(540,373)
\path(540,373)(990,1273)
\path(990,1273)(990,373)
\path(990,1273)(1440,373)
\path(1890,1273)(1890,373)
\path(1890,1273)(2340,373)
\dottedline{45}(990,373)(1890,1273)
\dottedline{45}(240,973)(390,673)
\blacken\path(309.502,766.915)(390.000,673.000)(363.167,793.748)(309.502,766.915)
\dottedline{45}(840,973)(690,673)
\blacken\path(716.833,793.748)(690.000,673.000)(770.498,766.915)(716.833,793.748)
\dottedline{45}(990,973)(990,673)
\blacken\path(960.000,793.000)(990.000,673.000)(1020.000,793.000)(960.000,793.000)
\dottedline{45}(1140,973)(1290,673)
\blacken\path(1209.502,766.915)(1290.000,673.000)(1263.167,793.748)(1209.502,766.915)
\dottedline{45}(1590,973)(1440,823)
\blacken\path(1503.640,929.066)(1440.000,823.000)(1546.066,886.640)(1503.640,929.066)
\dottedline{45}(1890,973)(1890,673)
\blacken\path(1860.000,793.000)(1890.000,673.000)(1920.000,793.000)(1860.000,793.000)
\dottedline{45}(2040,973)(2190,673)
\blacken\path(2109.502,766.915)(2190.000,673.000)(2163.167,793.748)(2109.502,766.915)
\path(990,1273)(1890,373)
\path(1590,673)(1740,523)
\blacken\path(1633.934,586.640)(1740.000,523.000)(1676.360,629.066)(1633.934,586.640)
\path(4065,1273)(4515,373)
\path(4215,973)(4365,673)
\blacken\path(4284.502,766.915)(4365.000,673.000)(4338.167,793.748)(4284.502,766.915)
\put(15,1423){\makebox(0,0)[lb]{\smash{{\SetFigFont{11}{13.2}{\familydefault}{\mddefault}{\updefault}$X_1$}}}}
\put(840,1423){\makebox(0,0)[lb]{\smash{{\SetFigFont{11}{13.2}{\familydefault}{\mddefault}{\updefault}$X_2$}}}}
\put(1740,1423){\makebox(0,0)[lb]{\smash{{\SetFigFont{11}{13.2}{\familydefault}{\mddefault}{\updefault}$X_3$}}}}
\put(390,73){\makebox(0,0)[lb]{\smash{{\SetFigFont{11}{13.2}{\familydefault}{\mddefault}{\updefault}$X_4$}}}}
\put(1740,73){\makebox(0,0)[lb]{\smash{{\SetFigFont{11}{13.2}{\familydefault}{\mddefault}{\updefault}$X_7$}}}}
\put(2265,73){\makebox(0,0)[lb]{\smash{{\SetFigFont{11}{13.2}{\familydefault}{\mddefault}{\updefault}$X_8$}}}}
\put(840,73){\makebox(0,0)[lb]{\smash{{\SetFigFont{11}{13.2}{\familydefault}{\mddefault}{\updefault}$X_5$}}}}
\put(1290,73){\makebox(0,0)[lb]{\smash{{\SetFigFont{11}{13.2}{\familydefault}{\mddefault}{\updefault}$X_6$}}}}
\put(3465,73){\makebox(0,0)[lb]{\smash{{\SetFigFont{11}{13.2}{\familydefault}{\mddefault}{\updefault}$X_4$}}}}
\put(4365,73){\makebox(0,0)[lb]{\smash{{\SetFigFont{11}{13.2}{\familydefault}{\mddefault}{\updefault}$X_6$}}}}
\put(4815,73){\makebox(0,0)[lb]{\smash{{\SetFigFont{11}{13.2}{\familydefault}{\mddefault}{\updefault}$X_7$}}}}
\put(5265,73){\makebox(0,0)[lb]{\smash{{\SetFigFont{11}{13.2}{\familydefault}{\mddefault}{\updefault}$X_8$}}}}
\put(2940,1423){\makebox(0,0)[lb]{\smash{{\SetFigFont{11}{13.2}{\familydefault}{\mddefault}{\updefault}$X_1$}}}}
\put(3690,1423){\makebox(0,0)[lb]{\smash{{\SetFigFont{11}{13.2}{\familydefault}{\mddefault}{\updefault}$(X_2 \setminus X_3) \cup (X_3 \setminus X_2)$}}}}
\end{picture}
}

%% file: treffig8-2.eepic
\setlength{\unitlength}{0.00083333in}
\begingroup\makeatletter\ifx\SetFigFont\undefined%
\gdef\SetFigFont#1#2#3#4#5{%
  \reset@font\fontsize{#1}{#2pt}%
  \fontfamily{#3}\fontseries{#4}\fontshape{#5}%
  \selectfont}%
\fi\endgroup%
{\renewcommand{\dashlinestretch}{30}
\begin{picture}(3963,1993)(0,-10)
\put(201,880){\ellipse{150}{150}}
\put(501,1180){\ellipse{150}{150}}
\put(801,1480){\ellipse{150}{150}}
\put(1101,1780){\ellipse{150}{150}}
\put(1401,580){\ellipse{150}{150}}
\put(501,580){\blacken\ellipse{150}{150}}
\put(501,580){\ellipse{150}{150}}
\put(801,880){\blacken\ellipse{150}{150}}
\put(801,880){\ellipse{150}{150}}
\put(1101,1180){\blacken\ellipse{150}{150}}
\put(1101,1180){\ellipse{150}{150}}
\put(1401,1480){\blacken\ellipse{150}{150}}
\put(1401,1480){\ellipse{150}{150}}
\put(201,1780){\blacken\ellipse{150}{150}}
\put(201,1780){\ellipse{150}{150}}
\put(2301,580){\ellipse{150}{150}}
\put(2601,280){\ellipse{150}{150}}
\put(3201,1780){\ellipse{150}{150}}
\put(3501,880){\ellipse{150}{150}}
\put(3801,1180){\ellipse{150}{150}}
\put(2301,1780){\blacken\ellipse{150}{150}}
\put(2301,1780){\ellipse{150}{150}}
\put(2601,880){\blacken\ellipse{150}{150}}
\put(2601,880){\ellipse{150}{150}}
\put(2901,580){\blacken\ellipse{150}{150}}
\put(2901,580){\ellipse{150}{150}}
\put(3501,1480){\blacken\ellipse{150}{150}}
\put(3501,1480){\ellipse{150}{150}}
\put(3801,280){\blacken\ellipse{150}{150}}
\put(3801,280){\ellipse{150}{150}}
\put(2901,1480){\ellipse{150}{150}}
\put(3201,1180){\blacken\ellipse{150}{150}}
\put(3201,1180){\ellipse{150}{150}}
\dottedline{45}(51,1930)(51,430)
\dottedline{45}(351,1930)(351,430)
\dottedline{45}(651,1930)(651,430)
\dottedline{45}(951,1930)(951,430)
\dottedline{45}(1251,1930)(1251,430)
\dottedline{45}(1551,1930)(1551,430)
\dottedline{45}(51,1930)(1551,1930)
\dottedline{45}(51,1630)(1551,1630)
\dottedline{45}(51,1330)(1551,1330)
\dottedline{45}(51,1030)(1551,1030)
\dottedline{45}(51,730)(1551,730)
\dottedline{45}(51,430)(1551,430)
\dottedline{45}(2151,1930)(2151,130)
\dottedline{45}(2151,1330)(3951,1330)
\dottedline{45}(2151,430)(3951,430)
\dottedline{45}(2451,1930)(2451,130)
\dottedline{45}(2751,1930)(2751,130)
\dottedline{45}(3051,1930)(3051,130)
\dottedline{45}(3351,1930)(3351,130)
\dottedline{45}(3651,1930)(3651,130)
\dottedline{45}(3951,1930)(3951,130)
\thicklines
\path(201,955)(201,1705)
\path(201,955)(201,1705)
\path(501,1105)(501,655)
\path(501,1105)(501,655)
\path(801,955)(801,1405)
\path(801,955)(801,1405)
\path(1101,1255)(1101,1705)
\path(1101,1255)(1101,1705)
\path(1401,1405)(1401,655)
\path(1401,1405)(1401,655)
\path(576,580)(1326,580)
\path(576,580)(1326,580)
\path(276,880)(426,880)
\path(276,880)(426,880)
\path(576,880)(726,880)
\path(576,880)(726,880)
\path(576,1180)(726,1180)
\path(576,1180)(726,1180)
\path(876,1180)(1026,1180)
\path(876,1180)(1026,1180)
\path(276,1780)(1026,1780)
\path(276,1780)(1026,1780)
\path(876,1480)(1026,1480)
\path(876,1480)(1026,1480)
\path(1176,1480)(1326,1480)
\path(1176,1480)(1326,1480)
\thinlines
\dottedline{45}(51,1930)(51,430)
\thicklines
\path(2301,1705)(2301,655)
\path(2301,1705)(2301,655)
\path(2601,805)(2601,355)
\path(2601,805)(2601,355)
\path(3501,1405)(3501,955)
\path(3501,1405)(3501,955)
\path(3801,1105)(3801,355)
\path(3801,1105)(3801,355)
\path(2376,1780)(3126,1780)
\path(2376,1780)(3126,1780)
\path(2676,880)(2826,880)
\path(2676,880)(2826,880)
\path(2976,880)(3426,880)
\path(2976,880)(3426,880)
\path(2376,580)(2526,580)
\path(2376,580)(2526,580)
\path(2676,580)(2826,580)
\path(2676,580)(2826,580)
\path(2676,280)(3726,280)
\path(2676,280)(3726,280)
\path(2901,655)(2901,1405)
\path(2901,655)(2901,1405)
\path(3201,1705)(3201,1255)
\path(3201,1705)(3201,1255)
\path(3276,1180)(3426,1180)
\path(3276,1180)(3426,1180)
\path(3576,1180)(3726,1180)
\path(3576,1180)(3726,1180)
\path(2976,1480)(3126,1480)
\path(2976,1480)(3126,1480)
\path(3276,1480)(3426,1480)
\path(3276,1480)(3426,1480)
\thinlines
\dottedline{45}(2151,1630)(3951,1630)
\dottedline{45}(2151,1930)(3951,1930)
\dottedline{45}(2151,730)(3951,730)
\dottedline{45}(2151,130)(3951,130)
\dottedline{45}(2151,1030)(3951,1030)
\put(1515,655){\makebox(0,0)[lb]{\smash{{\SetFigFont{10}{12.0}{\familydefault}{\mddefault}{\updefault}$0$}}}}
\put(1515,955){\makebox(0,0)[lb]{\smash{{\SetFigFont{10}{12.0}{\familydefault}{\mddefault}{\updefault}$0$}}}}
\put(1515,358){\makebox(0,0)[lb]{\smash{{\SetFigFont{10}{12.0}{\familydefault}{\mddefault}{\updefault}$0$}}}}
\put(1515,1255){\makebox(0,0)[lb]{\smash{{\SetFigFont{10}{12.0}{\familydefault}{\mddefault}{\updefault}$0$}}}}
\put(1515,1555){\makebox(0,0)[lb]{\smash{{\SetFigFont{10}{12.0}{\familydefault}{\mddefault}{\updefault}$0$}}}}
\put(1515,1855){\makebox(0,0)[lb]{\smash{{\SetFigFont{10}{12.0}{\familydefault}{\mddefault}{\updefault}$0$}}}}
\put(1215,358){\makebox(0,0)[lb]{\smash{{\SetFigFont{10}{12.0}{\familydefault}{\mddefault}{\updefault}$0$}}}}
\put(1215,1555){\makebox(0,0)[lb]{\smash{{\SetFigFont{10}{12.0}{\familydefault}{\mddefault}{\updefault}$0$}}}}
\put(1215,1855){\makebox(0,0)[lb]{\smash{{\SetFigFont{10}{12.0}{\familydefault}{\mddefault}{\updefault}$0$}}}}
\put(915,358){\makebox(0,0)[lb]{\smash{{\SetFigFont{10}{12.0}{\familydefault}{\mddefault}{\updefault}$0$}}}}
\put(915,1255){\makebox(0,0)[lb]{\smash{{\SetFigFont{10}{12.0}{\familydefault}{\mddefault}{\updefault}$0$}}}}
\put(915,1555){\makebox(0,0)[lb]{\smash{{\SetFigFont{10}{12.0}{\familydefault}{\mddefault}{\updefault}$1$}}}}
\put(915,1855){\makebox(0,0)[lb]{\smash{{\SetFigFont{10}{12.0}{\familydefault}{\mddefault}{\updefault}$0$}}}}
\put(615,955){\makebox(0,0)[lb]{\smash{{\SetFigFont{10}{12.0}{\familydefault}{\mddefault}{\updefault}$0$}}}}
\put(615,358){\makebox(0,0)[lb]{\smash{{\SetFigFont{10}{12.0}{\familydefault}{\mddefault}{\updefault}$0$}}}}
\put(615,1255){\makebox(0,0)[lb]{\smash{{\SetFigFont{10}{12.0}{\familydefault}{\mddefault}{\updefault}$1$}}}}
\put(615,1555){\makebox(0,0)[lb]{\smash{{\SetFigFont{10}{12.0}{\familydefault}{\mddefault}{\updefault}$1$}}}}
\put(615,1855){\makebox(0,0)[lb]{\smash{{\SetFigFont{10}{12.0}{\familydefault}{\mddefault}{\updefault}$0$}}}}
\put(315,655){\makebox(0,0)[lb]{\smash{{\SetFigFont{10}{12.0}{\familydefault}{\mddefault}{\updefault}$0$}}}}
\put(315,955){\makebox(0,0)[lb]{\smash{{\SetFigFont{10}{12.0}{\familydefault}{\mddefault}{\updefault}$1$}}}}
\put(315,358){\makebox(0,0)[lb]{\smash{{\SetFigFont{10}{12.0}{\familydefault}{\mddefault}{\updefault}$0$}}}}
\put(315,1255){\makebox(0,0)[lb]{\smash{{\SetFigFont{10}{12.0}{\familydefault}{\mddefault}{\updefault}$1$}}}}
\put(315,1555){\makebox(0,0)[lb]{\smash{{\SetFigFont{10}{12.0}{\familydefault}{\mddefault}{\updefault}$1$}}}}
\put(315,1855){\makebox(0,0)[lb]{\smash{{\SetFigFont{10}{12.0}{\familydefault}{\mddefault}{\updefault}$0$}}}}
\put(15,655){\makebox(0,0)[lb]{\smash{{\SetFigFont{10}{12.0}{\familydefault}{\mddefault}{\updefault}$0$}}}}
\put(15,955){\makebox(0,0)[lb]{\smash{{\SetFigFont{10}{12.0}{\familydefault}{\mddefault}{\updefault}$0$}}}}
\put(15,358){\makebox(0,0)[lb]{\smash{{\SetFigFont{10}{12.0}{\familydefault}{\mddefault}{\updefault}$0$}}}}
\put(15,1255){\makebox(0,0)[lb]{\smash{{\SetFigFont{10}{12.0}{\familydefault}{\mddefault}{\updefault}$0$}}}}
\put(15,1555){\makebox(0,0)[lb]{\smash{{\SetFigFont{10}{12.0}{\familydefault}{\mddefault}{\updefault}$0$}}}}
\put(15,1855){\makebox(0,0)[lb]{\smash{{\SetFigFont{10}{12.0}{\familydefault}{\mddefault}{\updefault}$0$}}}}
\put(576,655){\makebox(0,0)[lb]{\smash{{\SetFigFont{10}{12.0}{\familydefault}{\mddefault}{\updefault}$-1$}}}}
\put(876,655){\makebox(0,0)[lb]{\smash{{\SetFigFont{10}{12.0}{\familydefault}{\mddefault}{\updefault}$-1$}}}}
\put(876,955){\makebox(0,0)[lb]{\smash{{\SetFigFont{10}{12.0}{\familydefault}{\mddefault}{\updefault}$-1$}}}}
\put(1176,1255){\makebox(0,0)[lb]{\smash{{\SetFigFont{10}{12.0}{\familydefault}{\mddefault}{\updefault}$-1$}}}}
\put(1176,955){\makebox(0,0)[lb]{\smash{{\SetFigFont{10}{12.0}{\familydefault}{\mddefault}{\updefault}$-1$}}}}
\put(1176,655){\makebox(0,0)[lb]{\smash{{\SetFigFont{10}{12.0}{\familydefault}{\mddefault}{\updefault}$-1$}}}}
\put(2115,1855){\makebox(0,0)[lb]{\smash{{\SetFigFont{10}{12.0}{\familydefault}{\mddefault}{\updefault}$0$}}}}
\put(2115,55){\makebox(0,0)[lb]{\smash{{\SetFigFont{10}{12.0}{\familydefault}{\mddefault}{\updefault}$0$}}}}
\put(2115,1555){\makebox(0,0)[lb]{\smash{{\SetFigFont{10}{12.0}{\familydefault}{\mddefault}{\updefault}$0$}}}}
\put(2115,1255){\makebox(0,0)[lb]{\smash{{\SetFigFont{10}{12.0}{\familydefault}{\mddefault}{\updefault}$0$}}}}
\put(2115,955){\makebox(0,0)[lb]{\smash{{\SetFigFont{10}{12.0}{\familydefault}{\mddefault}{\updefault}$0$}}}}
\put(2115,655){\makebox(0,0)[lb]{\smash{{\SetFigFont{10}{12.0}{\familydefault}{\mddefault}{\updefault}$0$}}}}
\put(2115,355){\makebox(0,0)[lb]{\smash{{\SetFigFont{10}{12.0}{\familydefault}{\mddefault}{\updefault}$0$}}}}
\put(2415,955){\makebox(0,0)[lb]{\smash{{\SetFigFont{10}{12.0}{\familydefault}{\mddefault}{\updefault}$1$}}}}
\put(2415,655){\makebox(0,0)[lb]{\smash{{\SetFigFont{10}{12.0}{\familydefault}{\mddefault}{\updefault}$1$}}}}
\put(2415,355){\makebox(0,0)[lb]{\smash{{\SetFigFont{10}{12.0}{\familydefault}{\mddefault}{\updefault}$0$}}}}
\put(2415,55){\makebox(0,0)[lb]{\smash{{\SetFigFont{10}{12.0}{\familydefault}{\mddefault}{\updefault}$0$}}}}
\put(2415,1255){\makebox(0,0)[lb]{\smash{{\SetFigFont{10}{12.0}{\familydefault}{\mddefault}{\updefault}$1$}}}}
\put(2415,1855){\makebox(0,0)[lb]{\smash{{\SetFigFont{10}{12.0}{\familydefault}{\mddefault}{\updefault}$0$}}}}
\put(2415,1555){\makebox(0,0)[lb]{\smash{{\SetFigFont{10}{12.0}{\familydefault}{\mddefault}{\updefault}$1$}}}}
\put(2715,1255){\makebox(0,0)[lb]{\smash{{\SetFigFont{10}{12.0}{\familydefault}{\mddefault}{\updefault}$1$}}}}
\put(2715,1855){\makebox(0,0)[lb]{\smash{{\SetFigFont{10}{12.0}{\familydefault}{\mddefault}{\updefault}$0$}}}}
\put(2715,1555){\makebox(0,0)[lb]{\smash{{\SetFigFont{10}{12.0}{\familydefault}{\mddefault}{\updefault}$1$}}}}
\put(3015,1255){\makebox(0,0)[lb]{\smash{{\SetFigFont{10}{12.0}{\familydefault}{\mddefault}{\updefault}$1$}}}}
\put(3015,1855){\makebox(0,0)[lb]{\smash{{\SetFigFont{10}{12.0}{\familydefault}{\mddefault}{\updefault}$0$}}}}
\put(3015,1555){\makebox(0,0)[lb]{\smash{{\SetFigFont{10}{12.0}{\familydefault}{\mddefault}{\updefault}$1$}}}}
\put(3255,1255){\makebox(0,0)[lb]{\smash{{\SetFigFont{10}{12.0}{\familydefault}{\mddefault}{\updefault}$-1$}}}}
\put(3315,1855){\makebox(0,0)[lb]{\smash{{\SetFigFont{10}{12.0}{\familydefault}{\mddefault}{\updefault}$0$}}}}
\put(3315,1555){\makebox(0,0)[lb]{\smash{{\SetFigFont{10}{12.0}{\familydefault}{\mddefault}{\updefault}$0$}}}}
\put(3615,1255){\makebox(0,0)[lb]{\smash{{\SetFigFont{10}{12.0}{\familydefault}{\mddefault}{\updefault}$0$}}}}
\put(3615,1855){\makebox(0,0)[lb]{\smash{{\SetFigFont{10}{12.0}{\familydefault}{\mddefault}{\updefault}$0$}}}}
\put(3615,1555){\makebox(0,0)[lb]{\smash{{\SetFigFont{10}{12.0}{\familydefault}{\mddefault}{\updefault}$0$}}}}
\put(3915,355){\makebox(0,0)[lb]{\smash{{\SetFigFont{10}{12.0}{\familydefault}{\mddefault}{\updefault}$0$}}}}
\put(3915,955){\makebox(0,0)[lb]{\smash{{\SetFigFont{10}{12.0}{\familydefault}{\mddefault}{\updefault}$0$}}}}
\put(3915,655){\makebox(0,0)[lb]{\smash{{\SetFigFont{10}{12.0}{\familydefault}{\mddefault}{\updefault}$0$}}}}
\put(3615,355){\makebox(0,0)[lb]{\smash{{\SetFigFont{10}{12.0}{\familydefault}{\mddefault}{\updefault}$1$}}}}
\put(3615,955){\makebox(0,0)[lb]{\smash{{\SetFigFont{10}{12.0}{\familydefault}{\mddefault}{\updefault}$1$}}}}
\put(3615,655){\makebox(0,0)[lb]{\smash{{\SetFigFont{10}{12.0}{\familydefault}{\mddefault}{\updefault}$1$}}}}
\put(3315,355){\makebox(0,0)[lb]{\smash{{\SetFigFont{10}{12.0}{\familydefault}{\mddefault}{\updefault}$1$}}}}
\put(3315,655){\makebox(0,0)[lb]{\smash{{\SetFigFont{10}{12.0}{\familydefault}{\mddefault}{\updefault}$1$}}}}
\put(3015,355){\makebox(0,0)[lb]{\smash{{\SetFigFont{10}{12.0}{\familydefault}{\mddefault}{\updefault}$1$}}}}
\put(3015,955){\makebox(0,0)[lb]{\smash{{\SetFigFont{10}{12.0}{\familydefault}{\mddefault}{\updefault}$1$}}}}
\put(3015,655){\makebox(0,0)[lb]{\smash{{\SetFigFont{10}{12.0}{\familydefault}{\mddefault}{\updefault}$1$}}}}
\put(2715,355){\makebox(0,0)[lb]{\smash{{\SetFigFont{10}{12.0}{\familydefault}{\mddefault}{\updefault}$1$}}}}
\put(2715,955){\makebox(0,0)[lb]{\smash{{\SetFigFont{10}{12.0}{\familydefault}{\mddefault}{\updefault}$1$}}}}
\put(2715,655){\makebox(0,0)[lb]{\smash{{\SetFigFont{10}{12.0}{\familydefault}{\mddefault}{\updefault}$2$}}}}
\put(3915,1255){\makebox(0,0)[lb]{\smash{{\SetFigFont{10}{12.0}{\familydefault}{\mddefault}{\updefault}$0$}}}}
\put(3915,1855){\makebox(0,0)[lb]{\smash{{\SetFigFont{10}{12.0}{\familydefault}{\mddefault}{\updefault}$0$}}}}
\put(3915,1555){\makebox(0,0)[lb]{\smash{{\SetFigFont{10}{12.0}{\familydefault}{\mddefault}{\updefault}$0$}}}}
\put(2715,55){\makebox(0,0)[lb]{\smash{{\SetFigFont{10}{12.0}{\familydefault}{\mddefault}{\updefault}$0$}}}}
\put(3015,55){\makebox(0,0)[lb]{\smash{{\SetFigFont{10}{12.0}{\familydefault}{\mddefault}{\updefault}$0$}}}}
\put(3315,55){\makebox(0,0)[lb]{\smash{{\SetFigFont{10}{12.0}{\familydefault}{\mddefault}{\updefault}$0$}}}}
\put(3615,55){\makebox(0,0)[lb]{\smash{{\SetFigFont{10}{12.0}{\familydefault}{\mddefault}{\updefault}$0$}}}}
\put(3915,55){\makebox(0,0)[lb]{\smash{{\SetFigFont{10}{12.0}{\familydefault}{\mddefault}{\updefault}$0$}}}}
\put(3315,955){\makebox(0,0)[lb]{\smash{{\SetFigFont{10}{12.0}{\familydefault}{\mddefault}{\updefault}$1$}}}}
\end{picture}
}

%% file: 10_154.eepic
\setlength{\unitlength}{0.00083333in}
\begingroup\makeatletter\ifx\SetFigFont\undefined%
\gdef\SetFigFont#1#2#3#4#5{%
  \reset@font\fontsize{#1}{#2pt}%
  \fontfamily{#3}\fontseries{#4}\fontshape{#5}%
  \selectfont}%
\fi\endgroup%
{\renewcommand{\dashlinestretch}{30}
\begin{picture}(6011,2652)(0,-10)
\path(1962,1212)(1962,12)
\path(1887,2037)(1737,1887)
\blacken\path(1800.640,1993.066)(1737.000,1887.000)(1843.066,1950.640)(1800.640,1993.066)
\path(1887,1737)(1737,1587)
\blacken\path(1800.640,1693.066)(1737.000,1587.000)(1843.066,1650.640)(1800.640,1693.066)
\path(1587,1437)(1437,1287)
\blacken\path(1500.640,1393.066)(1437.000,1287.000)(1543.066,1350.640)(1500.640,1393.066)
\path(1287,1137)(1137,987)
\blacken\path(1200.640,1093.066)(1137.000,987.000)(1243.066,1050.640)(1200.640,1093.066)
\path(987,837)(837,687)
\blacken\path(900.640,793.066)(837.000,687.000)(943.066,750.640)(900.640,793.066)
\path(387,537)(237,387)
\blacken\path(300.640,493.066)(237.000,387.000)(343.066,450.640)(300.640,493.066)
\path(12,1212)(2112,1212)
\path(1962,1212)(1962,2412)
\path(5112,2412)(5112,12)
\path(4737,1662)(4587,1362)
\blacken\path(4613.833,1482.748)(4587.000,1362.000)(4667.498,1455.915)(4613.833,1482.748)
\path(5262,1212)(3312,1212)
\put(1887,2037){\makebox(0,0)[lb]{\smash{{\SetFigFont{12}{14.4}{\familydefault}{\mddefault}{\updefault}$\ZZ_2^1$}}}}
\put(1887,1737){\makebox(0,0)[lb]{\smash{{\SetFigFont{11}{13.2}{\familydefault}{\mddefault}{\updefault}$\ZZ_2^1$}}}}
\put(1587,1437){\makebox(0,0)[lb]{\smash{{\SetFigFont{11}{13.2}{\familydefault}{\mddefault}{\updefault}$\ZZ_2^4$}}}}
\put(1287,1137){\makebox(0,0)[lb]{\smash{{\SetFigFont{11}{13.2}{\familydefault}{\mddefault}{\updefault}$\ZZ_2^7$}}}}
\put(987,837){\makebox(0,0)[lb]{\smash{{\SetFigFont{11}{13.2}{\familydefault}{\mddefault}{\updefault}$\ZZ_2^4$}}}}
\put(687,537){\makebox(0,0)[lb]{\smash{{\SetFigFont{11}{13.2}{\familydefault}{\mddefault}{\updefault}$\ZZ_2^1$}}}}
\put(387,537){\makebox(0,0)[lb]{\smash{{\SetFigFont{11}{13.2}{\familydefault}{\mddefault}{\updefault}$\ZZ_2^1$}}}}
\put(87,237){\makebox(0,0)[lb]{\smash{{\SetFigFont{11}{13.2}{\familydefault}{\mddefault}{\updefault}$\ZZ_2^1$}}}}
\put(1512,2487){\makebox(0,0)[lb]{\smash{{\SetFigFont{11}{13.2}{\familydefault}{\mddefault}{\updefault}Alexander}}}}
\put(1662,2037){\makebox(0,0)[lb]{\smash{{\SetFigFont{8}{9.6}{\familydefault}{\mddefault}{\updefault}$0$}}}}
\put(162,537){\makebox(0,0)[lb]{\smash{{\SetFigFont{8}{9.6}{\familydefault}{\mddefault}{\updefault}$1$}}}}
\put(1062,612){\makebox(0,0)[lb]{\smash{{\SetFigFont{8}{9.6}{\familydefault}{\mddefault}{\updefault}$1$}}}}
\put(1287,912){\makebox(0,0)[lb]{\smash{{\SetFigFont{8}{9.6}{\familydefault}{\mddefault}{\updefault}$3$}}}}
\put(1887,1512){\makebox(0,0)[lb]{\smash{{\SetFigFont{8}{9.6}{\familydefault}{\mddefault}{\updefault}$1$}}}}
\put(312,2112){\makebox(0,0)[lb]{\smash{{\SetFigFont{11}{13.2}{\familydefault}{\mddefault}{\updefault}$E^1$}}}}
\put(2262,1137){\makebox(0,0)[lb]{\smash{{\SetFigFont{11}{13.2}{\familydefault}{\mddefault}{\updefault}Maslov}}}}
\put(1587,1212){\makebox(0,0)[lb]{\smash{{\SetFigFont{8}{9.6}{\familydefault}{\mddefault}{\updefault}$3$}}}}
\put(1587,1737){\makebox(0,0)[lb]{\smash{{\SetFigFont{11}{13.2}{\familydefault}{\mddefault}{\updefault}$\ZZ_2^1$}}}}
\put(5037,2037){\makebox(0,0)[lb]{\smash{{\SetFigFont{12}{14.4}{\familydefault}{\mddefault}{\updefault}$\ZZ_2^1$}}}}
\put(4737,1737){\makebox(0,0)[lb]{\smash{{\SetFigFont{12}{14.4}{\familydefault}{\mddefault}{\updefault}$\ZZ_2^1$}}}}
\put(4437,1137){\makebox(0,0)[lb]{\smash{{\SetFigFont{12}{14.4}{\familydefault}{\mddefault}{\updefault}$\ZZ_2^1$}}}}
\put(4512,1512){\makebox(0,0)[lb]{\smash{{\SetFigFont{8}{9.6}{\familydefault}{\mddefault}{\updefault}$1$}}}}
\put(3462,2112){\makebox(0,0)[lb]{\smash{{\SetFigFont{11}{13.2}{\familydefault}{\mddefault}{\updefault}$E^2$}}}}
\put(4662,2487){\makebox(0,0)[lb]{\smash{{\SetFigFont{11}{13.2}{\familydefault}{\mddefault}{\updefault}Alexander}}}}
\put(5412,1137){\makebox(0,0)[lb]{\smash{{\SetFigFont{11}{13.2}{\familydefault}{\mddefault}{\updefault}Maslov}}}}
\end{picture}
}

%% file: kino.eepic
\setlength{\unitlength}{0.00083333in}
\begingroup\makeatletter\ifx\SetFigFont\undefined%
\gdef\SetFigFont#1#2#3#4#5{%
  \reset@font\fontsize{#1}{#2pt}%
  \fontfamily{#3}\fontseries{#4}\fontshape{#5}%
  \selectfont}%
\fi\endgroup%
{\renewcommand{\dashlinestretch}{30}
\begin{picture}(5874,2836)(0,-10)
\put(139,1410){\ellipse{126}{126}}
\put(394,2682){\ellipse{126}{126}}
\put(648,2427){\ellipse{126}{126}}
\put(1157,2173){\ellipse{126}{126}}
\put(903,1157){\ellipse{126}{126}}
\put(1411,139){\ellipse{126}{126}}
\put(1665,394){\ellipse{126}{126}}
\put(1919,1664){\ellipse{126}{126}}
\put(2174,1919){\ellipse{126}{126}}
\put(2428,648){\ellipse{126}{126}}
\put(2683,902){\ellipse{126}{126}}
\put(394,1664){\blacken\ellipse{126}{126}}
\put(394,1664){\ellipse{126}{126}}
\put(648,394){\blacken\ellipse{126}{126}}
\put(648,394){\ellipse{126}{126}}
\put(903,1919){\blacken\ellipse{126}{126}}
\put(903,1919){\ellipse{126}{126}}
\put(1157,2682){\blacken\ellipse{126}{126}}
\put(1157,2682){\ellipse{126}{126}}
\put(1411,648){\blacken\ellipse{126}{126}}
\put(1411,648){\ellipse{126}{126}}
\put(1665,1410){\blacken\ellipse{126}{126}}
\put(1665,1410){\ellipse{126}{126}}
\put(1919,2427){\blacken\ellipse{126}{126}}
\put(1919,2427){\ellipse{126}{126}}
\put(2174,902){\blacken\ellipse{126}{126}}
\put(2174,902){\ellipse{126}{126}}
\put(2428,1157){\blacken\ellipse{126}{126}}
\put(2428,1157){\ellipse{126}{126}}
\put(2683,2173){\blacken\ellipse{126}{126}}
\put(2683,2173){\ellipse{126}{126}}
\put(139,139){\blacken\ellipse{126}{126}}
\put(139,139){\ellipse{126}{126}}
\put(3191,2682){\ellipse{126}{126}}
\put(3446,2427){\ellipse{126}{126}}
\put(3700,902){\ellipse{126}{126}}
\put(3955,1157){\ellipse{126}{126}}
\put(4209,1410){\ellipse{126}{126}}
\put(4463,2173){\ellipse{126}{126}}
\put(4717,1664){\ellipse{126}{126}}
\put(4972,1919){\ellipse{126}{126}}
\put(5226,394){\ellipse{126}{126}}
\put(5480,648){\ellipse{126}{126}}
\put(5735,139){\ellipse{126}{126}}
\put(3191,1664){\blacken\ellipse{126}{126}}
\put(3191,1664){\ellipse{126}{126}}
\put(3446,394){\blacken\ellipse{126}{126}}
\put(3446,394){\ellipse{126}{126}}
\put(3700,1919){\blacken\ellipse{126}{126}}
\put(3700,1919){\ellipse{126}{126}}
\put(3955,139){\blacken\ellipse{126}{126}}
\put(3955,139){\ellipse{126}{126}}
\put(4209,902){\blacken\ellipse{126}{126}}
\put(4209,902){\ellipse{126}{126}}
\put(4463,2682){\blacken\ellipse{126}{126}}
\put(4463,2682){\ellipse{126}{126}}
\put(4717,2427){\blacken\ellipse{126}{126}}
\put(4717,2427){\ellipse{126}{126}}
\put(4972,648){\blacken\ellipse{126}{126}}
\put(4972,648){\ellipse{126}{126}}
\put(5226,1157){\blacken\ellipse{126}{126}}
\put(5226,1157){\ellipse{126}{126}}
\put(5480,2173){\blacken\ellipse{126}{126}}
\put(5480,2173){\ellipse{126}{126}}
\put(5735,1410){\blacken\ellipse{126}{126}}
\put(5735,1410){\ellipse{126}{126}}
\dottedline{45}(3064,2809)(5862,2809)(5862,12)
	(3064,12)(3064,2809)
\dottedline{45}(3064,2555)(5862,2555)
\dottedline{45}(3064,2300)(5862,2300)
\dottedline{45}(3064,2046)(5862,2046)
\dottedline{45}(3064,1792)(5862,1792)
\dottedline{45}(3064,1537)(5862,1537)
\dottedline{45}(3064,1284)(5862,1284)
\dottedline{45}(3064,1029)(5862,1029)
\dottedline{45}(3064,775)(5862,775)
\dottedline{45}(3064,521)(5862,521)
\dottedline{45}(3064,266)(5862,266)
\dottedline{45}(3064,2809)(5862,2809)(5862,12)
	(3064,12)(3064,2809)
\dottedline{45}(3319,2809)(3319,12)
\dottedline{45}(3573,2809)(3573,12)
\dottedline{45}(3828,2809)(3828,12)
\dottedline{45}(4082,2809)(4082,12)
\dottedline{45}(4590,2809)(4590,12)
\dottedline{45}(4844,2809)(4844,12)
\dottedline{45}(5099,2809)(5099,12)
\dottedline{45}(5353,2809)(5353,12)
\dottedline{45}(5608,2809)(5608,12)
\dottedline{45}(4336,2809)(4336,12)
\dottedline{45}(12,2809)(2810,2809)(2810,12)
	(12,12)(12,2809)
\dottedline{45}(12,2555)(2810,2555)
\thicklines
\path(648,2364)(648,457)
\path(648,2364)(648,457)
\path(394,2619)(394,1728)
\path(394,2619)(394,1728)
\path(139,202)(139,1347)
\path(139,202)(139,1347)
\path(903,1220)(903,1856)
\path(903,1220)(903,1856)
\path(1157,2237)(1157,2619)
\path(1157,2237)(1157,2619)
\path(1411,202)(1411,584)
\path(1411,202)(1411,584)
\path(1665,457)(1665,1347)
\path(1665,457)(1665,1347)
\path(1919,1728)(1919,2364)
\path(1919,1728)(1919,2364)
\path(2174,965)(2174,1856)
\path(2174,965)(2174,1856)
\path(2428,711)(2428,1093)
\path(2428,711)(2428,1093)
\path(2683,965)(2683,2110)
\path(2683,965)(2683,2110)
\path(202,1410)(584,1410)
\path(202,1410)(584,1410)
\path(711,1410)(838,1410)
\path(711,1410)(838,1410)
\path(966,1410)(1602,1410)
\path(966,1410)(1602,1410)
\path(457,1664)(584,1664)
\path(457,1664)(584,1664)
\path(711,1664)(838,1664)
\path(711,1664)(838,1664)
\path(966,1664)(1856,1664)
\path(966,1664)(1856,1664)
\path(966,1919)(1856,1919)
\path(966,1919)(1856,1919)
\path(1983,1919)(2111,1919)
\path(1983,1919)(2111,1919)
\path(711,2427)(1093,2427)
\path(711,2427)(1093,2427)
\path(1220,2427)(1856,2427)
\path(1220,2427)(1856,2427)
\path(457,2682)(1093,2682)
\path(457,2682)(1093,2682)
\path(1220,2173)(1856,2173)
\path(1220,2173)(1856,2173)
\path(1983,2173)(2619,2173)
\path(1983,2173)(2619,2173)
\path(202,139)(1347,139)
\path(202,139)(1347,139)
\path(711,394)(1347,394)
\path(711,394)(1347,394)
\path(1475,394)(1602,394)
\path(1475,394)(1602,394)
\path(2238,902)(2365,902)
\path(2238,902)(2365,902)
\path(2492,902)(2619,902)
\path(2492,902)(2619,902)
\path(966,1157)(1602,1157)
\path(966,1157)(1602,1157)
\path(1729,1157)(2111,1157)
\path(1729,1157)(2111,1157)
\path(2238,1157)(2365,1157)
\path(2238,1157)(2365,1157)
\path(1475,648)(1602,648)
\path(1475,648)(1602,648)
\path(1729,648)(2365,648)
\path(1729,648)(2365,648)
\thinlines
\dottedline{45}(12,2300)(2810,2300)
\dottedline{45}(12,2046)(2810,2046)
\dottedline{45}(12,1792)(2810,1792)
\dottedline{45}(12,1537)(2810,1537)
\dottedline{45}(12,1284)(2810,1284)
\dottedline{45}(12,1029)(2810,1029)
\dottedline{45}(12,775)(2810,775)
\dottedline{45}(12,521)(2810,521)
\dottedline{45}(12,266)(2810,266)
\dottedline{45}(12,2809)(2810,2809)(2810,12)
	(12,12)(12,2809)
\dottedline{45}(266,2809)(266,12)
\dottedline{45}(521,2809)(521,12)
\dottedline{45}(775,2809)(775,12)
\dottedline{45}(1030,2809)(1030,12)
\dottedline{45}(1284,2809)(1284,12)
\dottedline{45}(1538,2809)(1538,12)
\dottedline{45}(1792,2809)(1792,12)
\dottedline{45}(2047,2809)(2047,12)
\dottedline{45}(2301,2809)(2301,12)
\dottedline{45}(2555,2809)(2555,12)
\thicklines
\path(3191,1728)(3191,2619)
\path(3191,1728)(3191,2619)
\path(3446,2364)(3446,457)
\path(3446,2364)(3446,457)
\path(3700,965)(3700,1856)
\path(3700,965)(3700,1856)
\path(3955,1093)(3955,202)
\path(3955,1093)(3955,202)
\path(4209,965)(4209,1347)
\path(4209,965)(4209,1347)
\path(4463,2237)(4463,2619)
\path(4463,2237)(4463,2619)
\path(4717,2364)(4717,1728)
\path(4717,2364)(4717,1728)
\path(4972,1856)(4972,711)
\path(4972,1856)(4972,711)
\path(5226,457)(5226,1093)
\path(5226,457)(5226,1093)
\path(5480,2110)(5480,711)
\path(5480,2110)(5480,711)
\path(5735,1347)(5735,202)
\path(5735,1347)(5735,202)
\path(4018,139)(5672,139)
\path(4018,139)(5672,139)
\path(3509,394)(3891,394)
\path(3509,394)(3891,394)
\path(4018,394)(5163,394)
\path(4018,394)(5163,394)
\path(3763,902)(3891,902)
\path(3763,902)(3891,902)
\path(4018,902)(4145,902)
\path(4018,902)(4145,902)
\path(3255,1664)(3382,1664)
\path(3255,1664)(3382,1664)
\path(3509,1664)(3636,1664)
\path(3509,1664)(3636,1664)
\path(3763,1664)(4654,1664)
\path(3763,1664)(4654,1664)
\path(4272,1410)(4908,1410)
\path(4272,1410)(4908,1410)
\path(5036,1410)(5417,1410)
\path(5036,1410)(5417,1410)
\path(5544,1410)(5672,1410)
\path(5544,1410)(5672,1410)
\path(4018,1157)(4145,1157)
\path(4018,1157)(4145,1157)
\path(4272,1157)(4908,1157)
\path(4272,1157)(4908,1157)
\path(5036,1157)(5163,1157)
\path(5036,1157)(5163,1157)
\path(3763,1919)(4654,1919)
\path(3763,1919)(4654,1919)
\path(4781,1919)(4908,1919)
\path(4781,1919)(4908,1919)
\path(4527,2173)(4654,2173)
\path(4527,2173)(4654,2173)
\path(4781,2173)(5417,2173)
\path(4781,2173)(5417,2173)
\path(3255,2682)(4400,2682)
\path(3255,2682)(4400,2682)
\path(3509,2427)(4400,2427)
\path(3509,2427)(4400,2427)
\path(4527,2427)(4654,2427)
\path(4527,2427)(4654,2427)
\path(5036,648)(5163,648)
\path(5036,648)(5163,648)
\path(5290,648)(5417,648)
\path(5290,648)(5417,648)
\end{picture}
}

%% file: hfk2.bbl
\begin{thebibliography}{1}

\bibitem{BP}Y. Bae, C.Y. Park. \emph{An upper bound of arc index of links}. Math. Proc. Cambridge Philos. Soc. {\bf 129} (2000) 491-500.

\bibitem{C} P.~R.~Cromwell. \emph{Arc presentations of knots and links.} Proc. 1995 Warsaw Knot Theory Conference.

\bibitem{MOS} C. Manolescu, P. Ozsv\'ath, S. Sarkar. \emph{A combinatorial description of knot Floer homology.} math.GT/0607691.

\bibitem{N} I. J. Nutt. \emph{Braid Index of Satellite Links}. Ph.D. Thesis, Liverpool University (1995).

\bibitem{OSz4} P. Ozsv\'ath, Z. Szab\'o. \emph{Heegaard Floer homology and alternating knots.} Geom. Topol. {\bf 7} (2003) 225-254.

\bibitem{OSz1} P. Ozsv\'ath, Z. Szab\'o. \emph{Holomorphic disks and genus bounds.} Geom. Topol. {\bf 8} (2004) 311-334.

\bibitem{OSz5} P. Ozsv\'ath, Z. Szab\'o. \emph{Holomorphic disks and topological invariants for closed three-manifolds}. math.SG/0101206.

\bibitem{OSz3} P. Ozsv{\'a}th and Z. Szab{\'o}. \emph{Holomorphic disks and knot invariants}. math.GT/0209056.

\bibitem{OSz6} P. Ozsv\'ath, Z. Szab\'o. \emph{Knot Floer homology and the four-ball genus.} Geom. Topol.  {\bf 7} (2003) 615-639.

\bibitem{OSz2} P. Ozsv\'ath, Z. Szab\'o. \emph{Knot Floer homology, genus bounds, and mutation.} math.GT/0303225.

\bibitem{R} J.~A.~Rasmussen. \emph{Floer homology and knot complements.} Ph.D. thesis, Harvard University (2003).

\bibitem{MC} M. Culler.  \emph{Gridlink: a tool for knot theorists.} www.math.uic.edu/ $\sim$culler/gridlink/

\bibitem{MH} M. Hedden. \emph{Knot Floer homology of Whitehead doubles.} math.GT/0606094.

\bibitem{HO} M. Hedden, P. Ording. \emph{The Ozsv{\'a}th-Szab{\'o} and Rasmussen concordance invariants are not equal.} math.GT/0512348.

\end{thebibliography}
